\pdfoutput=1
\RequirePackage{ifpdf}
\ifpdf % We are running pdfTeX in pdf mode
\documentclass[pdftex]{sigma}
\else
\documentclass{sigma}
\fi

\numberwithin{equation}{section}

\DeclareMathOperator{\longg}{long}

\newtheorem{Theorem}{Theorem}[section]
\newtheorem{Lemma}[Theorem]{Lemma}
\newtheorem{Proposition}[Theorem]{Proposition}
{\theoremstyle{definition}
\newtheorem{Remarks}[Theorem]{Remarks}
\newtheorem{Remark}[Theorem]{Remark}
\newtheorem{observations}[Theorem]{Observations}
}

\begin{document}

\allowdisplaybreaks

\renewcommand{\PaperNumber}{050}

\FirstPageHeading

\ShortArticleName{Geometric Aspects of the Painlev\'e Equations ${\rm PIII(D_6)}$ and ${\rm PIII(D_7)}$}

\ArticleName{Geometric Aspects of the Painlev\'e Equations\\
$\boldsymbol{\rm PIII(D_6)}$ and $\boldsymbol{\rm PIII(D_7)}$}

\Author{Marius VAN DER PUT and Jaap TOP}
\AuthorNameForHeading{M.~van der Put and J.~Top}
\Address{Johann Bernoulli Institute, University of Groningen,\\
P.O.~Box~407, 9700 AK Groningen, The Netherlands}
\Email{\href{mailto:mvdput@math.rug.nl}{mvdput@math.rug.nl}, \href{mailto:j.top@rug.nl}{j.top@rug.nl}}
\URLaddress{\url{http://www.math.rug.nl/~top/}}

\ArticleDates{Received October 15, 2013, in f\/inal form April 10, 2014; Published online April 23, 2014}

\Abstract{The Riemann--Hilbert approach for the equations ${\rm PIII(D_6)}$ and ${\rm PIII(D_7)}$ is studied in detail,
involving moduli spaces for connections and monodromy data, Okamoto--Painlev\'e varieties, the Painlev\'e property,
special solutions and explicit B\"acklund transformations.}


\Classification{14D20; 14D22; 34M55}

\section{Introduction}

The aim of this paper is a~study of the Painlev\'e equations ${\rm PIII(D_6)}$ and ${\rm PIII(D_7)}$ by means of
isomonodromic families in a~moduli space of connections of rank two on the projective line.
This contrasts the work of Okamoto et al.
on these third Painlev\'e equations, where the Hamiltonians are the main tool.
However, there is a~close relation between the two points of view, since the moduli spaces turn out to be the
Okamoto--Painlev\'e varieties.
We refer only to a~few items of the extensive literature on Okamoto--Painlev\'e varieties.
More details on Stokes matrices and the analytic classif\/ication of singularities can be found in~\cite{vdP-Si}.

A {\it rough sketch} of this Riemann--Hilbert method is as follows (see for some background~\cite{vdP-Sa} and see for
details concerning PI, PII, PIV~\cite{vdP2,vdP1,vdP-T}).
The starting point is a~family $\bf S$ of dif\/ferential modules $M$ of dimension 2 over $\mathbb{C}(z)$ with prescribed
singularities at f\/ixed points of~$\mathbb{P}^1$.
The type of singularities gives rise to a~monodromy set~$\mathcal{R}$ built out of ordinary monodromy, Stokes matrices
and `links'.
The map ${\bf S}\rightarrow \mathcal{R}$ associates to each module $M\in {\bf S}$ its monodromy data in~$\mathcal{R}$.
The f\/ibers of ${\bf S}\rightarrow \mathcal{R}$ are parametrized by $T\cong \mathbb{C}^*$ and there results a~bijection
${\bf S}\rightarrow \mathcal{R}\times T$.
The set $\bf S$ has a~priori no structure of algebraic variety.
A moduli space~$\mathcal{M}$ over~$\mathbb{C}$, whose set of closed points consists of certain connections of rank two
on the projective line, is constructed such that $\bf S$ coincides with $\mathcal{M}(\mathbb{C})$.
There results an analytic Riemann--Hilbert morphism ${\rm RH}:\mathcal{M}\rightarrow \mathcal{R}$.
The f\/ibers of~$RH$ are the isomonodromic families which give rise to solutions of the corresponding Painlev\'e equation.
The extended Riemann--Hilbert morphism ${\rm RH}^+:\mathcal{M}\rightarrow \mathcal{R}\times T$ is an analytic isomorphism.
From these constructions the Painlev\'e property for the corresponding Painlev\'e equation follows and the moduli space~$\mathcal{M}$ is identif\/ied with an Okamoto--Painlev\'e space.
Special properties of solutions of the Painlev\'e equations, such as special solutions, B\"acklund transformations etc.,
are derived from special points of~$\mathcal{R}$ and the natural automorphisms of~$\bf S$.

The above sketch needs many subtle ref\/inements.
One has to construct a~geometric monodromy space $\tilde{\mathcal{R}}$ (depending on the parameters of the Painlev\'e
equation) which is a~geometric quotient of the monodromy data.
The `link' involves (multi)summation and in order to avoid the singular directions one has to replace $T$ by its
universal covering $\tilde{T}\cong \mathbb{C}$.
In the construction of $\mathcal{M}$, one represents a~dif\/ferential module $M\in \bf{S}$ by a~connection on a~f\/ixed
vector bundle of rank two on $\mathbb{P}^1$.
This works well for the case $(1,-,1/2)$, described in Section~\ref{Section1}.
In case $(1,-,1)$, described in Section~\ref{Section2}, this excludes a~certain set of reducible modules $M\in \bf{S}$.
The moduli space of connections $\mathcal{M}$ is replaced by the topological covering
$\tilde{\mathcal{M}}=\mathcal{M}\times_T\tilde{T}$.
Now the main result is that the extended Riemann--Hilbert morphism ${\rm RH}^+:\tilde{\mathcal{M}}\rightarrow
\tilde{\mathcal{R}}\times \tilde{T}$ is a~well def\/ined analytic isomorphism.

Each family of linear dif\/ferential modules and each Painlev\'e equation has its own story.
For the family $(1,-,1/2)$, corresponding to $\rm{PIII(D_7)}$, the computations of~$\mathcal{R}$, $\mathcal{M}$ and the
B\"acklund transformations present no problems.
For the resonant case (i.e., $\alpha =\pm 1$ and $\theta \in \mathbb{Z}$, see Sections~\ref{Section1.1},~\ref{Section1.2}
and~\ref{Section1.5} for notation and the statement) there are algebraic solutions of $\rm{PIII(D_7)}$.
The spaces $\mathcal{R}(\alpha)$ with $\alpha =\pm i$ have a~special point corresponding to trivial Stokes matrices.
This leads to a~special solution $q$ for $\rm{PIII(D_7)}$ and $\theta \in \frac{1}{2}+\mathbb{Z}$, which is
transcendental according to~\cite{Oh}.
According to~\cite{KO}, $q$ is a~univalent function of $t$ and is a~meromorphic at $t=0$.

For the general case (i.e., $\alpha \neq \beta^{\pm 1}$) of the family $(1,-,1)$, corresponding to $\rm{PIII(D_6)}$,
the computations of $\mathcal{R}$, $\mathcal{M}$ present no problems.
The formulas for the B\"acklund transformations, derived from the automorphisms of $\bf S$, have denominators.
These originate from the complicated cases $\alpha =\beta^{\pm 1}$ and/or $\alpha=\pm 1$ where reducible connections
and/or resonance occur.
Isomonodromy for reducible connections produces Riccati solutions and resonance is related to algebraic solutions.

\section[The family $(1,-,1/2)$ and ${\rm PIII(D_7)}$]{The family $\boldsymbol{(1,-,1/2)}$ and $\boldsymbol{\rm PIII(D_7)}$}
\label{Section1}

In this section the set $\bf S$ consists of the (equivalence classes of the) pairs $(M,t)$ of type $(1,-,1/2)$
(see~\cite{vdP-Sa} for the terminology), corresponding to the Painlev\'e equation ${\rm PIII (D_7)}$.
The dif\/ferential module $M$ is given by $\dim M=2$, the second exterior power of $M$ is trivial, $M$ has two singular
points 0 and $\infty$.
The Katz invariant $r(0)=1$ and the generalized eigenvalues at $0$ are normalized to $\pm \frac{t}{2}z^{-1}$ with $t\in
T=\mathbb{C}^*$.
The singular point $\infty$ has Katz invariant $r(\infty )=1/2$ and generalized eigenvalues $\pm z^{1/2}$.
Further $(M,t)$ is equivalent to $(M',t')$ if $M$ is isomorphic to $M'$ and $t=t'$.

The Riemann--Hilbert approach to ${\rm PIII(D_7)}$ in this section dif\/fers from~\cite{vdP-Sa} in several ways.
The choice $(1,-,1/2)$ is made to obtain the classical formula for the Painlev\'e equation.
Further we consider pairs $(M,t)$ rather than modules $M$.
This is needed in order to distinguish the two generalized eigenvalues at $z=0$ and to obtain a~good monodromy space~$\mathcal{R}$.
Finally, the def\/initions of the topological monodromy and the `link' need special attention.

\subsection[The construction of the monodromy space $\mathcal{R}\rightarrow \mathcal{P}$]{The construction
of the monodromy space $\boldsymbol{\mathcal{R}\rightarrow \mathcal{P}}$}
\label{Section1.1}

This is rather subtle and we provide here the details.
Given is some $(M,t)\in {\bf S}$ and we write $\delta_M$ for the dif\/ferential operator on $M$.
First we f\/ix an isomorphism $\phi:\Lambda^2M\rightarrow (\mathbb{C}(z),z\frac{d}{dz})$.

1.~$\mathbb{C}((z))\otimes M$ has a~basis $F_1$, $F_2$ such that the operator $\delta_M$
has the matrix $\begin{pmatrix}-\omega&0\\0&\omega\end{pmatrix}$
with $\omega =\frac{tz^{-1}+\theta}{2}$.
We require that $\phi (F_1\wedge F_2)=1$.
Then $F_1$, $F_2$ is unique up to a~transformation $(F_1,F_2)\mapsto (\lambda F_1,\lambda^{-1}F_2)$.
The solution space $V(0)$ at $z=0$ has basis $f_1=e^{-\frac{t}{2}z^{-1}}z^{\theta
/2}F_1,f_2=e^{\frac{t}{2}z^{-1}}z^{-\theta /2}F_2$.
The formal monodromy and the two Stokes matrices have on the basis $f_1$, $f_2$
the matrices $\begin{pmatrix}\alpha&0\\0&\alpha^{-1}\end{pmatrix}$,
$\begin{pmatrix}1&0\\c_1&1\end{pmatrix}$, $\begin{pmatrix}1&c_2\\0&1\end{pmatrix}$, where $\alpha =e^{\pi i \theta}$.
Their product (in this order) is the topological monodromy $\operatorname{top}_0$ at $z=0$.

2.~$\mathbb{C}((z^{-1}))\otimes M$ has a~basis $E_1$, $E_2$ such that the operator $\delta_M$ has
the matrix $\begin{pmatrix}\frac{1}{4}&1\\z&\frac{-1}{4}\end{pmatrix}$ with respect to this basis.
Since this dif\/ferential operator is irreducible, the basis $E_1$, $E_2$ is unique up to a~transformation $(E_1,E_2)\mapsto
(\mu E_1,\mu E_2)$ with $\mu \in \mathbb{C}^*$.
We require that $\phi (E_1\wedge E_2)=1$ and then $\mu \in \{1,-1\}$.
The solution space $V(\infty)$ at $z=\infty$ has basis
\begin{gather*}
e_1=\frac{1}{\sqrt{-2}}z^{-1/4}e^{2z^{1/2}}\left(E_1+\left(-\frac{1}{2}+z^{1/2}\right)E_2\right),
\\
e_2=\frac{1}{\sqrt{-2}}z^{-1/4}e^{-2z^{1/2}}\left(E_1+\left(-\frac{1}{2}-z^{1/2}\right)E_2\right).
\end{gather*}
The formal monodromy and the Stokes matrix have on the basis $e_1$, $e_2$ the matrices $\begin{pmatrix}0&-i\\-i&0\end{pmatrix}$,
$\begin{pmatrix}1&0\\e&1\end{pmatrix}$.
Their product (in this order) is the topological monodromy $\operatorname{top}_\infty$ at $z=\infty$.

3.~The link $L:V(0)\rightarrow V(\infty)$ is a~linear map obtained from multisummation at $z=0$, analytic continuation
along a~path from $0$ to $\infty$ and the inverse of multisummation at $z=\infty$.
The matrix $\begin{pmatrix}\ell_1&\ell_2\\ \ell_3&\ell_4\end{pmatrix}$ of~$L$ with respect to the bases $e_1$, $e_2$ and $f_1$, $f_2$ has
determinant~1, due to the isomorphism $\phi: \Lambda^2M\rightarrow (\mathbb{C}(z),z\frac{d}{dz})$ and the careful
choices of the bases.
The relation $\begin{pmatrix}\alpha&\alpha c_2\\ \frac{c_1}{\alpha}&\frac{1+c_1c_2}{\alpha}\end{pmatrix}
=\operatorname{top}_0=L^{-1}\circ \operatorname{top}_\infty \circ L$ yields $\alpha=-i\ell_{14}e+i\ell_{12}-i\ell_{34}\neq 0$ where $\ell_{ij}:=\ell_i\ell_j$.
One observes that $\alpha$, $c_1$, $c_2$ are determined by $L$ and $\operatorname{top}_\infty$.
Thus the af\/f\/ine space, given by the above data and relations, has coordinate ring
\begin{gather*}
\mathbb{C}\left[e,\ell_1,\ell_2,\ell_3,\ell_4,\frac{1}{-\ell_{14}e+\ell_{12}-\ell_{34}}\right]/(\ell_{14}-\ell_{23}-1).
\end{gather*}

4.~The group $G$ generated by the base changes $(e_1,e_2){\mapsto} (-e_1,-e_2)$
and $(f_1,f_2){\mapsto} (\lambda f_1,\lambda^{-1}f_2)$, acts on this af\/f\/ine space.
The monodromy space~$\mathcal{R}$ is the categorical quotient by $G$ and has coordinate ring
\begin{gather*}
\mathbb{C}\left[e,\ell_{12},\ell_{14},\ell_{23},\ell_{34},\frac{1}{ -\ell_{14}e+\ell
_{12}-\ell_{34}}\right]/(\ell_{14}-\ell_{23}-1, \ell_{12}\ell_{34}-\ell_{14}\ell_{23}).
\end{gather*}
This is in fact a~{\it geometric quotient}.
The morphism $\mathcal{R}\rightarrow \mathcal{P}=\mathbb{C}^*$ is given by $(e,\ell
_{12},\ell_{14},\ell_{23},\ell_{34})\!\mapsto$ $\alpha:=-i\ell_{14}e+i\ell_{12}-i\ell_{34}\neq 0$.
For a~suitable linear change of the variables, the f\/ibers $\mathcal{R}(\alpha)$ of $\mathcal{R}\rightarrow \mathcal{P}$
are nonsingular, af\/f\/ine cubic surfaces with three lines at inf\/inity, given by the equation
$x_1x_2x_3 +x_1^2+x_2^2+\alpha x_1 +x_2=0$.
A detailed calculation resulting in this equation is presented in~\cite[Section~3.5]{vdP-Sa}.

\begin{Lemma}%\label{L1.1}
The space $\mathcal{R}(\alpha)$ is simply connected.
\end{Lemma}

\begin{proof}
We remove from $\mathcal{R}(\alpha)$ the line $L:=\{(0,0,x_3)\,|\, x_3\in \mathbb{C}\}$ and project onto
$\mathbb{C}^2\setminus S$ by $(x_1,x_2,x_3)\mapsto (x_1,x_2) $.
Here $S$ is the union of $\{(0,x_2)\,|\, x_2\neq -1\}$ and $\{(x_1,0)\,|\, x_1\neq -\alpha\}$.
If $x_1x_2\neq 0$, then the f\/iber is one point.
If $x_1x_2=0$, then the f\/iber is an af\/f\/ine line.
Since $\mathbb{C}^2\setminus S$ is simply connected, $\mathcal{R}(\alpha)\setminus L$ is simply connected.
Then $\mathcal{R}(\alpha)$ is simply connected, too.
\end{proof}

{\it Remark on the differential Galois group}.
The dif\/ferential Galois group of a~module $M$, with $(M,t)\in {\bf S}$, can be considered as algebraic subgroup of ${\rm GL}(V(0))$.
It is the smallest algebraic subgroup containing the local dif\/ferential Galois group $G_0\subset {\rm GL}(V(0))$ at
$z=0$ and $L^{-1}G_\infty L$, where $G_\infty\subset {\rm GL}(V(\infty ))$ is the local dif\/ferential Galois group at
$z=\infty$.
Now $G_0$ is generated (as algebraic group) by the exponential torus
$\left\{\begin{pmatrix}s_1&0\\ 0&s_1^{-1}\end{pmatrix}\Bigg|s_1\in \mathbb{C}^*\right\}$,
the formal monodromy $\begin{pmatrix}\alpha&0\\ 0&\alpha^{-1}\end{pmatrix}$
and the Stokes maps $\begin{pmatrix}1&0\\ c_1&1\end{pmatrix}$, $\begin{pmatrix}1&c_2\\ 0&1\end{pmatrix}$.
The group $G_\infty$ is (as algebraic group) generated by the exponential torus
$\left\{\begin{pmatrix}s_2&0\\ 0&s_2^{-1}\end{pmatrix}\Bigg|s_2\in \mathbb{C}^*\right\}$, the formal monodromy $\begin{pmatrix}0&-i\\ -i
&0\end{pmatrix}$ and the Stokes map $\begin{pmatrix}1&0\\ e
&1\end{pmatrix}$.
This easily implies that the dif\/ferential Galois group is ${\rm SL}(2,\mathbb{C})$.
In particular, $M$ is irreducible and the same holds for the dif\/ferential module $\mathbb{C}(\sqrt[m]{z})\otimes M$ over
$\mathbb{C}(\sqrt[m]{z})$ for any $m\geq 2$.

{\it The construction needed to define the topological monodromy and the link}.

For the def\/inition of the link and the topological monodromies we have to choose nonsingular directions for the two
multisummations and a~path from $0$ to $\infty$.
At $z=\infty$ the singular direction does not depend on $t\in T$ and we can take a~f\/ixed nonsingular direction.
However, at $z=0$, the singular directions for $t\in T=\mathbb{C}^*$, $t=|t|e^{i\phi}$ are $\phi$ and $\pi +\phi$ and
they vary with $t$.
Thus we cannot use a~f\/ixed path from $0$ to $\infty$.
In order to obtain a~globally def\/ined map $L:V(0)\rightarrow V(\infty)$ we replace $T=\mathbb{C}^*$ by its universal
covering $\tilde{T}=\mathbb{C}\rightarrow T$, $\tilde{t}\mapsto e^{\tilde{t}}$.
The elements of $\tilde{T}\cong \mathbb{R}_{> 0}\times \mathbb{R}$ are written as $\tilde{t}=|t|e^{i\phi}$.
Consider the path $\tilde{z}=re^{id(r)}$, $0<r<\infty $, with
$d(r)=(\phi-\frac{\pi}{2})\frac{1}{1+r}+\frac{\pi}{2}\frac{r}{1+r}$ on the universal covering of $\mathbb{P}^1 \setminus
\{0,\infty\}$.
Now $L$ is def\/ined by summation at $z=0$ in the direction $\phi -\frac{\pi}{2}$, followed by analytic continuation along
the above path and f\/inally the inverse of the summation at $z=\infty$ in the direction $\frac{\pi}{2}$.

Write $\tilde{\bf S}={\bf S}\times_T\tilde{T}$.
The elements of $\tilde{\bf S}$ are the pairs $(M,\tilde{t})$ with $(M,e^{\tilde{t}})\in {\bf S}$.
For the elements in $\tilde{\bf S}$ the link and the monodromy at $z=0$ are def\/ined as above.
Since~$\mathcal{R}$ is a~geometric quotient,~\cite[Theorem~1.9]{vdP-Sa} implies:

{\it The above map $\tilde{\bf S}\rightarrow \mathcal{R}\times \tilde{T}$ is bijective}.

Fix $\alpha \in \mathcal{P}=\mathbb{C}^*$.
Let ${\bf S}(\alpha)$, $\tilde{{\bf S}}(\alpha)$ be the subsets of $\bf S$ and $\tilde{\bf S}$, consisting of the pairs
$(M,t)$ and $(M,\tilde{t})$ which have $\begin{pmatrix}\alpha&0\\0&\frac{1}{\alpha}\end{pmatrix}$ as formal monodromy at $z=0$.

{\it The map $\tilde{\bf S}(\alpha)\rightarrow \mathcal{R}(\alpha )\times \tilde{T}$ is bijective}, since
$\tilde{\bf S}\rightarrow \mathcal{R}\times \tilde{T}$ is bijective.

\subsection[The construction of the moduli space $\mathcal{M}(\theta)$]
{The construction of the moduli space $\boldsymbol{\mathcal{M}(\theta)}$}
\label{Section1.2}

Fix $\theta$ with $e^{\pi i \theta}=\alpha$.
The moduli space $\mathcal{M}(\theta )$ is obtained by replacing each $(M,t)\in {\bf S}(\alpha )$ by a~certain
connection $(\mathcal{V},\nabla )$ on $\mathbb{P}^1$.
This connection is uniquely determined by the data: Its generic f\/iber is $M$; $\nabla_{z\frac{d}{dz}}$ is formally
equivalent at $z=0$ to $z\frac{d}{dz}+\begin{pmatrix}\omega&0\\ 0&-\omega\end{pmatrix}$
with $\omega =\frac{tz^{-1}+\theta}{2}$ and is formally equivalent at $z=\infty$ to
$z\frac{d}{dz}+\begin{pmatrix}-\frac{3}{4}&1\\z&-\frac{1}{4}\end{pmatrix}$.
It follows that $\Lambda^2(\mathcal{V},\nabla)$ is isomorphic to $(O(-1),d)$.
Since $(\mathcal{V},\nabla)$ is irreducible one has that $\mathcal{V}\cong O\oplus O(-1)$ and the vector bundle
$\mathcal{V}$ is identif\/ied with $Oe_1+O(-[\infty])e_2$.

Then $D:=\nabla_{z\frac{d}{dz}}:\mathcal{V}\rightarrow O([0]+[\infty ])\otimes \mathcal{V}$ has with respect to
$e_1$, $e_2$ the matrix $\begin{pmatrix}a&b\\ c
&-a\end{pmatrix}$ with $a=a_{-1}z^{-1}+a_0+a_1z$, $c=c_{-1}z^{-1}+c_0$ and $b=b_{-1}z^{-1}+b_0+b_1z+b_2z^2$.

{\it The condition at $z=0$} is $a^2+bc\in \big(\frac{tz^{-1}+\theta}{2}\big)^2+\mathbb{C}[[z]]$, equivalently
\begin{gather*}
a_{-1}^2+b_{-1}c_{-1}=\frac{t^2}{4},\qquad 2a_{-1}a_0+b_{-1}c_0+b_0c_{-1}=\frac{t\theta}{2}.
\end{gather*}

{\it The condition at $z=\infty$} is $a^2+a+bc=z+\mathbb{C}[[z^{-1}]]$, equivalently
\begin{gather*}
a_1^2+b_2c_0=0,\qquad 2a_1a_0+a_1+b_2c_{-1}+b_1c_0=1.
\end{gather*}

The space, given by the above variables and relations has to be divided by the action of the group $\{e_1\mapsto
e_1$, $e_2\mapsto \lambda e_2+(x_0+x_1z)e_1\}$ (with $\lambda \in \mathbb{C}^*$, $x_0,x_1 \in \mathbb{C}$) of automorphisms
of the vector bundle.
Using the standard forms below one sees that this is a~good geometric quotient.

A {\it standard form} for $c_{-1}\neq 0$ is $z\frac{d}{dz}+\begin{pmatrix}a_1z&b\\z^{-1}+c_0&-a_1z\end{pmatrix}$
with $b=b_{-1}z^{-1}+\dots +b_2z^2$ and equations
\begin{gather*}
a_1^2+b_2c_0=0,\qquad a_1+b_2+b_1c_0=1,\qquad b_{-1}=\frac{t^2}{4},\qquad \frac{t^2}{4}c_0+b_0=\frac{t\theta}{2}.
\end{gather*}

A {\it standard form} for $c_0\neq 0$ is $z\frac{d}{dz}+\begin{pmatrix}a_{-1}z^{-1}&b\\ c_{-1}z^{-1}+1&-a_{-1}z^{-1}\end{pmatrix}$ with equations
\begin{gather*}
b_2=0,\qquad b_1=1,\qquad a_{-1}^2+b_{-1}c_{-1}=\frac{t^2}{4},\qquad b_{-1}+b_0c_{-1}=\frac{t\theta}{2}.
\end{gather*}
By gluing the two standard forms, one obtains the nonsingular moduli space $\mathcal{M}(\theta)$.
The map $\mathcal{M}(\theta )\rightarrow {\bf S}(\alpha)$, where $\alpha =e^{i\theta}$, is a~bijection.

{\bf Observation.}
After scaling some variables one sees that $\mathcal{M}(\theta)$ is the union of two open af\/f\/ine spaces $U_1\times T$
and $U_2\times T$, where $U_1$ is given by the variables $a_1$, $b_1$, $c_0$ and the relation $a_1^2+(1-a_1-b_1c_0)c_0=0$, and
$U_2$ is given by the variables $a_{-1}$, $b_0$, $c_{-1}$ and the relation
$a_{-1}^2+(\frac{\theta}{2}-b_0c_{-1})c_{-1}-\frac{1}{4}=0$.

Let $U_{12}\subset U_1$ be def\/ined by $c_0\neq 0$ and $U_{21}\subset U_2$ by $c_{-1}\neq 0$.
The gluing of $U_1\times T$ and $U_2\times T$ is def\/ined by the isomorphism $U_{12}\times T\rightarrow U_{21}\times T$
obtained by suitable base changes in the group $\{e_1\mapsto e_1$, $e_2\mapsto \lambda e_2+(x_0+x_1z)e_1\}$.

Using the two projections $U_1\rightarrow \mathbb{C}$, $(a_1,b_1,c_0)\mapsto c_0$ and $U_2\rightarrow \mathbb{C}$, $
(a_{-1},b_0,c_{-1})\mapsto c_{-1}$ one f\/inds that $U_1$ and $U_2$ are simply connected.

{\it Define $M(\theta ):=f^{-1}(1)$, where $f:\mathcal{M}(\theta )\rightarrow T$ is the canonical morphism}.
The space $M(\theta)$ is simply connected since it is the union of the two simply connected spaces $U_1$, $U_2$.

{\it The universal covering of $\mathcal{M}(\theta )$ is $\tilde{\mathcal{M}}(\theta)=\mathcal{M}(\theta)\times
_T\tilde{T}$}.
Indeed, it is the union of the two simply connected spaces $U_1\times \tilde{T}$ and $U_2\times \tilde{T}$.
Using the explicitly def\/ined link $L$ one obtains a~globally def\/ined analytic morphism $\tilde{\mathcal{M}}(\theta
)\rightarrow \mathcal{R}(\alpha)\times \tilde{T}$ which is bijective (and thus an analytic isomorphism, see~\cite{KK}).
Indeed, $\mathcal{M}(\theta)\rightarrow \bf{S}(\alpha)$ and $\tilde{\bf{S}}(\alpha)\rightarrow \mathcal{R}(\alpha)\times
\tilde{T}$ are bijections.

\begin{Theorem}
\label{T1.2}
Let $\theta \in \mathbb{C}$, $\alpha =e^{\pi i \theta}$.
The extended Riemann--Hilbert map $\tilde{\mathcal{M}}( \theta )\rightarrow \mathcal{R}(\alpha)\times \tilde{T}$, with
$\tilde{T}\cong \mathbb{C}$, is a~well defined analytic isomorphism.
\end{Theorem}

{\bf Comment.} The existence of an analytic isomorphism as in Theorem~\ref{T1.2} is called the ``geometric
Painlev\'e property'' in~\cite{SaitoInaba}.
They prove this property for a~number of Painlev\'e equations under a~restriction on the parameters (loc.~cit.,
Theorem~6.3).
We prove it here for ${\rm PIII(D_7)}$ and in Sections~\ref{section2.3.2} and~\ref{section2.4.4} below for ${\rm PIII(D_6)}$ without any restriction.

\subsection{Isomonodromy and the Okamoto-Painlev\'e space}

The calculation is done on the `chart' $c_0\neq 0$ and $q$ is supposed to be invertible.
The data for the operator $z\frac{d}{dz}+A$ are
\begin{gather*}
c_{-1}=-q,
\qquad
b_{-1}=q^{-1}\left(a_{-1}^2-\frac{t^2}{4}\right),
\qquad
b_0=-\frac{t\theta}{2}q^{-1}+q^{-2}\left(a_{-1}^2-\frac{t^2}{4}\right),
\\
b_1=1,
\qquad
b_2=0.
\end{gather*}
Now $q$ and $a_{-1}$ are functions of $t$ and the family is isomonodromic if there is an operator $\frac{d}{dt}+B$
commuting with $z\frac{d}{dz} +A$.
Equivalently $A'=\overset{\rm o}{B}+[A,B]$, where $'$ denotes $\frac{d}{dt}$ and $\overset{\rm o}{}$ denotes $z\frac{d}{dz}$.

The Lie algebra $\frak{sl}_2$ has standard basis $H=\begin{pmatrix}1&0\\ 0&-1\end{pmatrix}$,
$E_1=\begin{pmatrix}0&1\\ 0&0\end{pmatrix}$, $E_2=\begin{pmatrix}0&0 \\ 1&0\end{pmatrix}$.
One writes $A=a_{-1}z^{-1}H+bE_1+(-qz^{-1}+1)E_2$ and $B=B_HH+B_1E_1+B_2E_2$ with $B_*=\sum\limits_{i=-1}^2B_{*,i}z^i$ for $*=H,1,2$
and $B_{*,i}$ only depending on $t$.
Using the Lie algebra structure one obtains the equations:
\begin{alignat*}{3}
& a_{-1}'z^{-1}=\overset{\rm o}{B_H}+B_2(z+b_0+b_{-1}z^{-1})-B_1(-qz^{-1}+1),
\qquad &&  (H)
\\
& b_0'+b_{-1}'z^{-1}=\overset{\rm o}{B_1}+ 2B_1a_{-1}z^{-1}-2B_H(z+b_0+b_{1}z^{-1}),
\qquad &&  (E_1)
\\
& -q'z^{-1}=\overset{\rm o}{B_2}-2B_2a_{-1}z^{-1}+2B_H(-qz^{-1}+1).
\qquad &&  (E_2)
\end{alignat*}
By Maple one obtains the system
\begin{gather*}
q'=\frac{q+2a_{-1}}{t},\qquad a_{-1}'=\frac{-t^2-\theta t q+4a_{-1}^2+2qa_{-1}+2q^3}{2tq}
\end{gather*}
and f\/inally
\begin{gather*}
q''=\frac{(q')^2}{q}-\frac{q'}{t}-\frac{\theta}{t}+\frac{2q^2}{t^2}-\frac{1}{q}.
\end{gather*}

We note that the change $q=-Q$, $t=-T$ brings this equation in the form
\begin{gather*}
\overset{**}{Q}=\frac{(\overset{*}{Q})^2}{Q}-\frac{\overset{*}{Q}}{T}-\frac{\theta}{T}-\frac{2Q^2}{T^2}-\frac{1}{Q}\qquad \text{
with notation}  \quad \overset{*}{-}=\frac{d-}{dT}.
\end{gather*}
This is the standard form for ${\rm PIII'(D_7)}$ (see~\cite{OKSO}).
As in~\cite{vdP2,vdP1,vdP-T} one obtains:

\begin{Theorem}
\label{T1.3}
The equation ${\rm PIII(D_7)}$ has the Painlev\'e property.
The analytic fibration $\tilde{t}:\tilde{\mathcal{M}}(\theta )\rightarrow \tilde{T}=\mathbb{C}$ with its foliation
$\{{\rm RH}^{-1}(r)\,|\, r\in \mathcal{R}(\alpha)\}$, where ${\rm RH}:\tilde{\mathcal{M}}(\theta )\rightarrow \mathcal{R}(\alpha)$ is
the Riemann--Hilbert map, is isomorphic to the Okamoto--Painlev\'e space for the equation ${\rm PIII(D_7)}$ with
parameter $\theta$.

Moreover, $M(\theta)\cong \mathcal{R}(\alpha)$ is the space of initial values.
\end{Theorem}

\subsection[Automorphisms of ${\bf S}$ and B\"acklund transforms]{Automorphisms of ${\bf S}$ and B\"acklund transforms}

The automorphism $s_1$ of $\bf S$ is def\/ined by $s_1(M,t)=(M,-t)$.
The induced action on~$\mathcal{R}$ leaves all data invariant except for interchanging the basis vectors $f_1$, $f_2$ of
$V(0)$.
As a~consequence $\alpha$ is mapped to $\alpha^{-1}$.

The automorphism $s_2$ of $\bf S$ is def\/ined by $s_2(M,t)=(N\otimes M,-t)$.
Here $N=\mathbb{C}(z)b$ is the dif\/ferential module given by $\delta (b)=\frac{1}{2}b$.
Starting with a~local presentation $z\frac{d}{dz}+\begin{pmatrix}\omega&0\\ 0&-\omega\end{pmatrix}$
with $\omega =\frac{ tz^{-1}+\theta}{2}$ of $M$ at $z=0$ one obtains,
after conjugation with $\begin{pmatrix}z&0\\0&1\end{pmatrix}$,
the local presentation $z\frac{d}{dz}+\begin{pmatrix}-\tau&0\\0&\tau\end{pmatrix}$,
with $\tau=\frac{tz^{-1}-\theta +1}{2}$, of $N\otimes M$ at $z=0$.

Starting with a~local presentation $D:=z\frac{d}{dz}+\begin{pmatrix}\frac{1}{4}&1\\ z
&-\frac{1}{4}\end{pmatrix}$ of $M$ at $z=\infty$ (say on the basis $E_1$, $E_2$, described in Section~\ref{Section1.1},  part~2), one obtains the
local presentation $z\frac{d}{dz}+\begin{pmatrix}\frac{3}{4}&1\\ z
&\frac{1}{4}\end{pmatrix}$ of $N\otimes M$ at $z=\infty$.
This is the matrix of $D$ with respect to the basis $zE_2$, $E_1$.
The induced action of $s_2$ on~$\mathcal{R}$ maps $\alpha$ to $-\alpha^{-1}$, the formal monodromy at $\infty$ is
multiplied by $-1$ and the Stokes data are essentially unchanged.

The group of automorphisms of ${\bf S}$, generated by $s_1$, $s_2$, has order 4.
The B\"acklund transformations are the lifts of the elements of this group to isomorphisms (preserving isomonodromy)
between various moduli spaces $\tilde{\mathcal{M}}( \theta )$.

$s_1^+:\tilde{\mathcal{M}}(\theta )\rightarrow \tilde{\mathcal{M}}(- \theta )$ is the obvious lift of $s_1$,
given by $\tilde{t}\mapsto \tilde{t}+\pi i$, $\theta \mapsto -\theta$.
Further any solution $q(\tilde{t})$ of ${\rm PIII(D_7)}$ for the parameter $\theta$ is mapped to the solution
$q(\tilde{t}+\pi i)$ for the parameter $-\theta$.

$s_2^+:\tilde{\mathcal{M}}(\theta )\rightarrow \tilde{\mathcal{M}}(1-\theta )$ is the obvious lift of $s_2$
with $\tilde{t}\mapsto \tilde{t}+\pi i$.
The formula for $s_2^+$ is not obvious and its computation is given below.

Put $B:=(s_1^+)^2=(s_2^+)^2$.
Then $B$ is the automorphism of $\tilde{\mathcal{M}(\theta)}=M(\theta)\times \tilde{T}$, which is the identity on
$M(\theta)$ and $B:\tilde{t}\mapsto \tilde{t}+2\pi i$.

The group $\langle s_1^+,s_2^+\rangle $ generated by $s_1^+$, $s_2^+$ (for their action on $\theta$, $\tilde{t}$) has $\langle B\rangle \cong \mathbb{Z}$
as normal subgroup and $ \langle s_1^+,s_2^+\rangle /\langle B\rangle $ is the af\/f\/ine Weyl group of type~$A_1$.

{\it Computation of the B\"acklund transformation $s_2^+$.}
A point $\xi \in \mathcal{M}(\theta )$, lying in the af\/f\/ine open subset def\/ined by $c_0\neq 0$ and $c_1\neq 0$, is
represented by the operator in standard form $z\frac{d}{dz}+\begin{pmatrix}az^{-1}&b\\ 1-qz^{-1}&-az^{-1}\end{pmatrix}$, where
$b=z-\frac{t\theta}{2q}+\frac{a^2-\frac{t^2}{4}}{q^2}+\frac{a^2-\frac{t^2}{4}}{q}z^{-1}$.
The map $s_2^+$ changes this operator into $z\frac{d}{dz}+A$, where $A$ is obtained from the above matrix by $t\mapsto
-t$ and adding $\begin{pmatrix}\frac{1}{2}&0\\ 0&\frac{1}{2}\end{pmatrix}$.
The point $s_2^+(\xi)\in \mathcal{M}(1-\theta)$ is supposed to be represented by the operator $z\frac{d}{dz}+\tilde{A}$,
where $\tilde{A}=\begin{pmatrix}\tilde{a}z^{-1}&\tilde{b}\\ 1-\tilde{q}z^{-1}&-\tilde{a}z^{-1}\end{pmatrix}$ with
$\tilde{b}=z-\frac{t(1-\theta)}{2\tilde{q}}+\frac{\tilde{a}^2-\frac{t^2}{4}}{\tilde{q}^2}+\frac{\tilde{a}^2-\frac{t^2}{4}}{\tilde{q}}z^{-1}$.
Since the two matrix dif\/ferential operators represent the same irreducible dif\/ferential module over $\mathbb{C}(z)$,
there is a~$T\in {\rm GL}(2,\mathbb{C}(z))\neq 0$, unique up to multiplication by a~constant, such that
$(z\frac{d}{dz}+A)T=T(z\frac{d}{dz}+\tilde{A})$.
A local computation shows that $T$ has the form $T_{-2}z^{-2}+T_{-1}z^{-1}+T_0\neq 0$ with constant matrices $T_*$.
The $\tilde{a}$, $\tilde{q}$ and the entries of the $T_*$ are the unknows in the identity
$(z\frac{d}{dz}+A)(T_{-2}z^{-2}+T_{-1}z^{-1}+T_0)=(T_{-2}z^{-2}+T_{-1}z^{-1}+T_0)(z\frac{d}{dz}+\tilde{A})$.
A Maple computation yields
\begin{gather*}
\tilde{q}=-\frac{t(\theta q+2a-t)}{2q^2}
\end{gather*}
and
\begin{gather*}
\tilde{a}=\frac{t(4a^2-4at+2aq+2\theta a~q+q^2\theta +t^2-tq\theta -qt-2q^3 )}{4q^3}.
\end{gather*}
The isomorphism $s_2^+$ respects the foliations.
For a~leaf one has $q'=\frac{q+2a}{t}$ and substitution in the f\/irst formula produces $\tilde{q}=-\frac{t(\theta
t+tq'-q-t)}{2q^2}$ for this B\"acklund transformation on solutions of ${\rm PIII(D_7)}$.

The B\"acklund transformation $s_2^+s_1^+$ maps a~solution $q$ for the parameter $\theta$ to the solution
$\frac{t(\theta q-2a+t)}{2q^2}$, with $a=\frac{tq'-q}{2}$, with parameter $1+\theta$.

\subsection{Remarks}\label{Section1.5}

1.~One considers for $(M,t)\in {\bf S}$ the connection $(\mathcal{V}_0,\nabla )$ with generic f\/ibre $M$ and the local data
$z\frac{d}{dz}+\begin{pmatrix}\frac{1}{4}&z\\ 1&-\frac{1}{4}\end{pmatrix}$
at $z=\infty$ and $z\frac{d}{dz}+\begin{pmatrix}\omega&0 \\ 0&-\omega\end{pmatrix}$ with $\omega =\frac{tz^{-1}+\theta}{2}$ at $z=0$.
The second exterior power of $(\mathcal{V}_0,\nabla )$ is $(O,d)$ and thus $\mathcal{V}_0$ has degree 0.
Since $(\mathcal{V}_0,\nabla )$ is irreducible, there are two possibilities for $\mathcal{V}_0$ namely $O\oplus O$ and
$O(1)\oplus O(-1)$.
Suppose that $\mathcal{V}_0\cong O(1)\oplus O(-1)$.
Then one can identify $\mathcal{V}_0$ with $O([\infty])e_1+O(-[\infty])e_2$ and for a~good choice of $e_1$, $e_2$ one
obtains the operator $\nabla_{z\frac{d}{dz}}= z\frac{d}{dz}+\begin{pmatrix}0&b\\z^{-1}&0\end{pmatrix}$ with $b=z^{2}+\cdots$.
One concludes that the locus of the modules $M$ with $\mathcal{V}_0=O(1)\oplus O(-1)$ is the set of the closed points of
$q^{-1}(\infty )$ for the map $q:=-\frac{c_0}{c_1}: \mathcal{M}(\theta )\rightarrow \mathbb{P}^1$ (see also
Section~\ref{Section2.3.4}).

2.~{\it Algebraic solutions of} ${\rm PIII(D_7)}$.
One easily f\/inds the algebraic solution(s) $q$ with $q^3=\frac{t^2}{2}$ for ${\rm PIII(D_7)}$ with $\theta =0$.
Using the B\"acklund transformations one f\/inds an algebraic solution for ${\rm PIII(D_7)}$ for every $\theta \in
\mathbb{Z}$.
According to~\cite{Oh, OKSO}, these are all the algebraic solutions of ${\rm PIII(D_7)}$.

More precisely, $q_j(\tilde{t})=e^{2\pi ij/3}\frac{e^{2\tilde{t}/3}}{\sqrt[3]{2}}$, $j=0,1,2$ are algebraic solution
for $\theta=0$.
We note that $q_1(\tilde{t})=q_0(\tilde{t}+4\pi i)$ and $q_2(\tilde{t})=q_0(\tilde{t}+2\pi i)$.
Since $\alpha=1$ and $\operatorname{top}^3=1$, these solutions are mapped to a~single point of $\mathcal{R}(1)$ corresponding to
$c_1c_2=-3$, $e=-i$ and certain values for the invariants $\ell_{12}$, $\ell_{14}$, $\ell_{23}$, $\ell_{34}$ (which we cannot make explicit).
The isomonodromic family for this solution $q$ is
\begin{gather*}
z\frac{d}{dz}+\begin{pmatrix}-\frac{q}{6}z^{-1}&b\\-qz^{-1}+1&\frac{q}{6}z^{-1}\end{pmatrix}
\qquad
\text{with}
\quad
b=z+\left(\frac{1}{36}-\frac{q}{2}\right)+q\left(\frac{1}{36}-\frac{q}{2}\right)z^{-1}
\end{gather*}
and
\begin{gather*}
q^3=\frac{t^2}{2}.
\end{gather*}
It is not clear what makes this family and the corresponding point of $\mathcal{R}(1)$ so special.

3.~{\it Special solutions of $\rm{PIII(D_7)}$}.
Consider an isomonodromy family for which the Stokes matrices are trivial, i.e., $c_1=c_2=e=0$.
Then $\alpha =i$ or $\alpha =-i$.
In the f\/irst case one computes that $\ell_{12}=\ell_{14}=\frac{1}{2}$, $\ell_{23}=\ell_{34}=-\frac{1}{2}$ and one f\/inds
a~unique point of $\mathcal{R}(i)$ and a~special solu\-tion~$q(\tilde{t})$ of $\rm{PIII(D_7)}$ for $\theta =\frac{1}{2}$.
Using B\"acklund transformations one obtains a~similar special solution for any $\theta \in \frac{1}{2}+\mathbb{Z}$.
Y.~Ohyama informed us that the condition $c_1=c_2=0$ implies that the corresponding solution $q$ of ${\rm PIII(D_7)}$ is
a~univalent function of $t$ and is meromorphic at $t=0$.
Further $\theta \in \frac{1}{2}+\mathbb{Z}$ is equivalent to $e=0$.
See~\cite{KO} for details.

\section[The family $(1,-,1)$]{The family $\boldsymbol{(1,-,1)}$}\label{Section2}

\subsection{Definition of the family}

The set $\bf S$ consists of the equivalence classes of pairs $(M,t)$, where $M$ is a~dif\/ferential module $M$ over
$\mathbb{C}(z)$ and $t\in \mathbb{C}^*$ such that: $\dim M=2$, $\Lambda^2M$ is the trivial module, $M$ has two
singulari\-ties~$0$ and~$\infty$, both singularities have Katz invariant $1$, the (generalized) eigenvalues are normalized
to $\pm \frac{t}{2}z^{-1}$ at $0$ and $\pm \frac{t}{2}z$ at $\infty$.
Further, two pairs $(M_1,t_1)$ and $(M_2,t_2)$ are called equivalent if there exists an isomorphism $M_1\rightarrow M_2$
and $t_1=t_2$.

As in Section~\ref{Section1},
we will have to replace $T$ by its universal covering $\tilde{T}=\mathbb{C}\rightarrow T$, $\tilde{t}\mapsto e^{\tilde{t}}$.
Write $\tilde{\bf S}=\bf{S}\times_T\tilde{T}$.
Def\/ine for $\alpha, \beta \in \mathbb{C}^*$ the subset $\bf{S}(\alpha, \beta)$ of $\bf{S}$ consisting of the pairs
$(M,t)$ such that $\mathbb{C}((z))\otimes M$ is represented by $z\frac{d}{dz}+\begin{pmatrix}\omega&0\\ 0&-\omega\end{pmatrix}$ with
$\omega = \frac{tz^{-1}+\theta_0}{2}$, $\alpha=e^{\pi i \theta_0}$ and $\mathbb{C}((z^{-1}))\otimes M$ is represented by
$z\frac{d}{dz}+\begin{pmatrix}\tau &0\\ 0&-\tau\end{pmatrix}$ with $\tau=\frac{tz+\theta_\infty}{2}$, $\beta=e^{\pi i \theta_\infty}$.
Further $\tilde{\bf S}(\alpha, \beta )=\bf{S}(\alpha, \beta )\times_T \tilde{T}$.

\subsection{The monodromy space}\label{Section2.2}

For $(M,t)\in {\bf S}$, the {\it monodromy data} are given by (compare~\cite{vdP-Sa}): the symbolic solutions spaces~$V(0)$ and~$V(\infty)$ at $z=0$ and $z=\infty$ (including formal monodromies and Stokes matrices) and the link
$L:V(0)\rightarrow V(\infty)$.
We make this more explicit.

The module $\mathbb{C}((z))\otimes M$ has a~basis $E_1$, $E_2$ with $\delta (E_1\wedge E_2)=0$ and $\delta
E_1=-\frac{tz^{-1}+\theta_0}{2}E_1$, $\delta E_2=\frac{tz^{-1}+\theta_0}{2}E_2$.
We note that $t$ is used to distinguish between $E_1$ and $E_2$.
This basis is unique up to a~transformation $E_1\mapsto c_1z^mE_1$, $E_2\mapsto c_2z^{-m}E_2$ with $c_1,c_2\in\mathbb{C}^*$,
$m\in \mathbb{Z}$.
After f\/ixing $\theta_0$, the $E_1$, $E_2$ are unique up to multiplication by constants.
The symbolic solution space $V(0)$ at $z=0$ is $\mathbb{C}e_1+\mathbb{C}e_2$, with
$e_1=e^{-\frac{t}{2}z^{-1}+\frac{\theta_0}{2}\log z}E_1$ and $e_2=e^{+\frac{t}{2}z^{-1}-\frac{\theta_0}{2}\log z}E_2$.

Now $\alpha =e^{\pi i\theta_0}$ is well def\/ined and does not depend on the choices for $E_1$, $E_2$.

Similarly, $\mathbb{C}((z^{-1}))\otimes M=\mathbb{C}((z^{-1})F_1\oplus \mathbb{C}((z^{-1}))F_2$ with $\delta (F_1\wedge
F_2)=0$, $\delta F_1=-\frac{tz+\theta_\infty}{2}F_1$ and $\delta F_2=\frac{tz+\theta_\infty}{2}F_2$.
The space $V(\infty)$ has basis $f_1=e^{\frac{tz}{2}+\frac{\theta_\infty}{2}\log z}F_1$ and
$f_2=e^{-\frac{tz}{2}-\frac{\theta_\infty}{2}\log z}F_2$ over $\mathbb{C}$.
Moreover, $\beta =e^{\pi i \theta_\infty}$.

For the basis $e_1$, $e_2$ of $V(0)$, the formal monodromy and the Stokes matrices are:
\begin{gather*}
\begin{pmatrix}\alpha&0\\ 0&\frac{1}{\alpha}\end{pmatrix},
\qquad
\begin{pmatrix}1&0\\ a_1&1\end{pmatrix},
\qquad
\begin{pmatrix}1&a_2\\0&1\end{pmatrix}
\qquad
\text{with product}
\qquad
\begin{pmatrix}\alpha&\alpha a_2\\ \frac{a_1}{\alpha}&\frac{1+a_1a_2}{\alpha}\end{pmatrix}.
\end{gather*}
This product is the topological monodromy $\operatorname{top}_0$ at $z=0$.

For the basis $f_1$, $f_2$ of $V(\infty)$, the formal monodromy and the Stokes matrices are:
\begin{gather*}
\begin{pmatrix}\beta&0\\ 0&\frac{1}{\beta}\end{pmatrix},
\qquad
\begin{pmatrix}1&0\\ b_1&1\end{pmatrix},
\qquad
\begin{pmatrix}1&b_2\\ 0&1\end{pmatrix}
\qquad
\text{with product}
\qquad
\begin{pmatrix}\beta&\beta b_2\\ \frac{b_1}{\beta}&\frac{1+b_1b_2}{\beta}\end{pmatrix}.
\end{gather*}
This is the topological monodromy $\operatorname{top}_\infty$ at $z=\infty$.

The link $L:V(0)\rightarrow V(\infty)$ with matrix $\begin{pmatrix}\ell_1&\ell_2\\ \ell_3&\ell_4\end{pmatrix}$ has determinant~1.

The relations are given by the matrix equality
\begin{gather*}
\begin{pmatrix}\beta&\beta b_2\\ \frac{b_1}{\beta}&\frac{1+b_1b_2}{\beta}\end{pmatrix}
=L\circ\begin{pmatrix}\alpha&\alpha a_2\\ \frac{a_1}{\alpha}&\frac{1+a_1a_2}{\alpha}\end{pmatrix}\circ L^{-1}.
\end{gather*}
In particular, $\beta =\ell_1\ell_4\alpha +\ell_2\ell_4\frac{a_1}{\alpha}-\ell_1\ell_3\alpha a_2-\ell_2\ell
_3\frac{1+a_1a_2}{\alpha}$.
This def\/ines a~variety $\mathcal{T}$, given by the variables $\alpha, a_1,a_2,\ell_1,\dots,\ell_4$
with the only restrictions $\ell_1\ell_4-\ell_2\ell_3=1$, $\alpha \neq 0$ and $\ell_1\ell_4\alpha +\ell_2\ell
_4\frac{a_1}{\alpha}-\ell_1\ell_3\alpha a_2-\ell_2\ell_3\frac{1+a_1a_2}{\alpha}\neq 0$.

For f\/ixed values of $\alpha,\beta \in \mathbb{C}^*$
we obtains a~variety $\mathcal{T}(\alpha,\beta)$ def\/ined by the
variables $a_1,a_2,\ell_1$, $\dots,\ell_4$ and the relations:
\begin{gather*}
\ell_1\ell_4-\ell_2\ell_3=1 \qquad \text{and}\qquad \beta =\ell_1\ell_4\alpha +\ell_2\ell_4\frac{a_1}{\alpha}-\ell_1\ell_3\alpha
a_2-\ell_2\ell_3\frac{1+a_1a_2}{\alpha}.
\end{gather*}

The group $\mathbb{G}_m\times \mathbb{G}_m$ acts on $\mathcal{T}$ and $\mathcal{T}(\alpha, \beta )$,
by base change $(e_1,e_2,f_1,f_2)\mapsto (\gamma e_1, \gamma^{-1} e_2,\delta f_1$, $\delta^{-1}f_2)$.
The categorical quotient of $\mathcal{T}$ by $\mathbb{G}_m\times \mathbb{G}_m$ is $\mathcal{R}\rightarrow \mathcal{P}$
with parameter space $\mathcal{P}=\mathbb{C}^*\times \mathbb{C}^*$ given by $(\alpha,\beta)$.
This is a~family of af\/f\/ine cubic surfaces $\mathcal{R}(\alpha,\beta)$ (this is the categorical quotient of
$\mathcal{T}(\alpha,\beta )$) given by the equation
\begin{gather*}
x_1x_2x_3+x_1^2+x_2^2+(1+\alpha \beta)x_1+(\alpha +\beta )x_2+\alpha \beta =0,
\end{gather*}
where
\begin{gather*}
x_1=\ell_1\ell_4-1,
\qquad
x_2=\alpha a_2\ell_1\ell_3-\alpha \ell_1\ell_4,
\qquad
x_3=\frac{1+a_1a_2}{\alpha}+\alpha.
\end{gather*}

{\it Observation:} $\mathcal{R}(\alpha,\beta)$ is simply connected if $\alpha \neq \beta^{\pm 1}$.

Def\/ine $U$ by removing the two lines $\{(0,-\beta,*)\}$, $\{(-1,0,*)\}$ from $\mathcal{R}(\alpha,\beta)$.
The image of the projection $(x_1,x_2,x_3)\in U\mapsto (x_1,x_2)\in \mathbb{C}^2$ is $\mathbb{C}^*\times
\mathbb{C}^*\cup \{ (0,-\alpha), (-\alpha \beta,0)\}$.
This image is simply connected.
For $x_1x_2\neq 0$, the f\/iber is one point.
For $x_1x_2=0$, the f\/iber is an af\/f\/ine line.
It follows that $U$ is simply connected and thus $\mathcal{R}(\alpha,\beta)$ is simply connected, too.

{\it Definition of the link and the topological monodromies}.
A~construction similar to the one in Section~\ref{Section1} is needed for the def\/inition of the link.
For $\tilde{t}=|t|e^{i\phi}\in \tilde{T}$, $\phi \in \mathbb{R}$ the singular directions at $z=0$ and $z=\infty$ are
$\phi$, $\phi -\pi$ and $-\phi$, $\pi-\phi$.
On the universal covering of $\mathbb{P}^1\setminus \{0,\infty\}$ one considers the path
$\tilde{z}=re^{id(r)}$, $0<r<\infty$ with $d(r)=\frac{1}{1+r}(\phi -\frac{\pi}{2})+\frac{r}{1+r}(\frac{\pi}{2}-\phi)$.
The link $L:V(0)\rightarrow V(\infty)$ is def\/ined by (multi)summation at zero in the direction $\phi -\frac{\pi}{2}$,
followed by analytic continuation along the above path and f\/inally the inverse of (multi)summation in the direction
$\frac{\pi}{2}-\phi$ at inf\/inity.
Now the map $\tilde{\bf S}(\alpha,\beta )\rightarrow \mathcal{R}(\alpha,\beta)\times \tilde{T}$ is well def\/ined.

{\it In the general case, i.e., $\alpha \neq \beta^{\pm 1}$}, the space $\mathcal{R}(\alpha,\beta)$ is the
geometric quotient of $\mathcal{T}(\alpha,\beta)$ and this space is nonsingular.
Therefore the natural map $\tilde{\bf S}(\alpha,\beta)\rightarrow \mathcal{R}(\alpha,\beta)\times \tilde{T}$ is
a~bijection.

Let $(M,\tilde{t})\in \tilde{{\bf S}}(\alpha, \beta)$.
Then $M$ is reducible if and only its monodromy data (in $\mathcal{R}(\alpha,\beta)$) is reducible.
Further $\mathcal{R}(\alpha,\beta)$ contains reducible monodromy data if and only if $\alpha=\beta^{\pm 1}$.
Thus ${\bf S}(\alpha,\beta)$ contains reducible modules if and only if $\alpha =\beta^{\pm 1}$.
We will f\/irst investigate the general case $\alpha \neq \beta^{\pm 1}$.
The special case $\alpha =\beta^{\pm 1}$ presents many dif\/f\/iculties and will be handled later on.

\subsection[The general case $\alpha \neq \beta^{\pm 1}$]{The general case $\boldsymbol{\alpha\neq\beta^{\pm 1}}$}

\subsubsection[The moduli space $\mathcal{M}(\theta_0,\theta_\infty)$]{The moduli space $\boldsymbol{\mathcal{M}(\theta_0,\theta_\infty)}$}
\label{Section2.3.1}

Fix $\theta_0$, $\theta_\infty$ with $\alpha =e^{i\pi \theta_0}$, $\beta =e^{i\pi \theta_\infty}$.
We will construct a~moduli space $\mathcal{M}(\theta_0,\theta_\infty)$ of connections
on the bundle $Oe_1\oplus O(-[0])e_2$ on $\mathbb{P}^1$ such that the map, which associates to a~connection in this space its generic f\/iber,
is a~bijection $\mathcal{M}(\theta_0,\theta_\infty)\rightarrow {\bf S}(\alpha,\beta)$.

The elements of the {\it set} $\mathcal{M}(\theta_0,\theta_\infty)$ are the connections $\nabla:\mathcal{V}\rightarrow
\Omega (2[0]+2[\infty])\otimes \mathcal{V}$ def\/ined by: the generic f\/iber $M$ satisf\/ies $(M,t)\in {\bf
S}(\alpha,\beta)$; the invariant lattices at $z=0$ and $z=\infty$ are given by the local matrix dif\/ferential operators
$z\frac{d}{dz}+\left(
\begin{matrix}
\frac{tz^{-1}+\theta_0}{2}&0
\\
0&-\frac{tz^{-1}+\theta_0}{2}+1
\end{matrix}
\right)$
and
$z\frac{d}{dz}+\left(
\begin{matrix}
\frac{tz+\theta_\infty}{2}&0
\\
0&-\frac{tz +\theta_\infty}{2}
\end{matrix}
\right)$.
Since the second exterior power of $M$ is trivial, the degree of $\mathcal{V}$ is $-1$.
By assumption $M$ is irreducible and therefore the type of $\mathcal{V}$ is $O\oplus O(-1)$ and one can identify~$\mathcal{V}$ with $Oe_1\oplus O(-[0])e_2$.
By construction the map $\mathcal{M}(\theta_0,\theta_\infty)\rightarrow {\bf S}(\alpha, \beta)$ is bijective.

The operator $\nabla_{z\frac{d}{dz}}$ has, with respect to the basis $\{e_1,e_2\}$,
the form $z\frac{d}{dz}+\begin{pmatrix}a&b\\ c&-a\end{pmatrix}$ with $a=a_{-1}z^{-1}+a_0+a_1z$, $b=b_{-2}z^{-2}+\dots +b_1z$, $c=c_0+c_1z$.
The lattice condition at $z=0$ is equivalent to $a(a-1)+bc\in(\frac{tz^{-1}+\theta_0}{2})^2-(\frac{tz^{-1}+\theta_0}{2})+\mathbb{C}[[z]]$.
This leads to the equations
\begin{gather*}
a_{-1}^2+b_{-2}c_0=\frac{t^2}{4},
\qquad
2a_{-1}a_0-a_{-1}+b_{-2}c_1+b_{-1}c_0=t\left(\frac{\theta_0}{2}-\frac{1}{2}\right).
\end{gather*}
The lattice condition at $z=\infty$ is equivalent to $a^2+bc\in (\frac{tz+\theta_\infty}{2})^2+\mathbb{C}[[z^{-1}]]$
and one f\/inds the equations
\begin{gather*}
a_1^2+b_1c_1=\frac{t^2}{4},
\qquad
2a_0a_1+b_0c_1+b_1c_0=\frac{t\theta_\infty}{2}.
\end{gather*}

Def\/ine the (quasi-)af\/f\/ine space $\mathcal{A}(\theta_0,\theta_\infty)$ by the variables $a_*$, $b_*$, $c_*$, $t$, the four
equations and the open condition $(c_0,c_1)\neq (0,0)$.
For the general case $\alpha \neq \beta^{\pm 1}$, the module $M$ and the corresponding connection are irreducible and
thus $(c_0,c_1)\neq (0,0)$ holds (see Remarks~\ref{R2.2}.3).

The space $\mathcal{A}(\theta_0,\theta_\infty)$ is divided out by the group of the automorphisms of the bundle
$Oe_1\oplus O(-[0])e_2$.
This action amounts to dividing $\mathcal{A}(\theta_0,\theta_\infty)$ by the group $G$ of transformations $e_1\mapsto
e_1$, $e_2\mapsto \lambda e_2+ (\gamma z^{-1}+\delta )e_1$.

\begin{Proposition}%\label{P2.1}
The quotient of $\mathcal{A}(\theta_0,\theta_\infty)$ by $G$ is geometric and has no singularities.
\end{Proposition}
\begin{proof}
The open subset $\mathcal{A}(\theta_0,\theta_\infty)_1$, def\/ined by $c_1\neq 0$, contains the closed `standard subset'
$ST_1$, given by the connections
\begin{gather*}
\begin{split}
& z\frac{d}{dz}+\left(
\begin{matrix}
a_{-1}z^{-1}&b
\\
z+c_0&-a_{-1}z^{-1}
\end{matrix}
\right)
\qquad
\text{with}
\qquad
b=b_{-2}z^{-2}+\dots +b_1z,
\\
& b_1=\frac{t^2}{4},\qquad b_0=-\frac{t^2}{4}c_0+\frac{t\theta_\infty}{2},\qquad b_{-2}=a_{-1}-b_{-1}c_0+t\left(\frac{\theta_0}{2}-\frac{1}{2}\right).
\end{split}
\end{gather*}
and the equation
\begin{gather*}
a_{-1}^2+c_0\left(a_{-1}-b_{-1}c_0+t\left(\frac{\theta_0}{2}-\frac{1}{2}\right)\right)-\frac{t^2}{4}=0.
\end{gather*}
The natural morphism $G\times ST_1\rightarrow \mathcal{M}(\theta_0,\theta_\infty)_1$ is an isomorphism.

Similarly, let $ST_0$ denote the closed subset of the open subset $\mathcal{A}(\theta_0,\theta_\infty)_0$, def\/ined by
$c_0\neq 0$, be given by the connections
\begin{gather*}
z\frac{d}{dz}+\left(
\begin{matrix}
a_{1}z&b
\\
c_1z+1&-a_1z
\end{matrix}
\right)
\qquad
\text{with}
\qquad
b=b_{-2}z^{-2}+\dots +b_1z,
\\
b_{-2}=\frac{t^2}{4},\qquad b_{-1}=t\left(\frac{\theta_0}{2}-\frac{1}{2}\right)-\frac{t^2}{4}c_1,\qquad b_1=\frac{t\theta_\infty}{2}-b_0c_1
\end{gather*}
and the equation
\begin{gather*}
a_1^2+\left(\frac{t\theta_\infty}{2}-b_0c_1\right)c_1-\frac{t^2}{4}=0.
\end{gather*}
Again $G\times ST_0\rightarrow \mathcal{A}(\theta_0,\theta_\infty)_0$ is an isomorphism.
The quotient of $\mathcal{A}(\theta_0,\theta_\infty)$ by $G$ is obtained by gluing the two non singular spaces $ST_1$
and $ST_0$.
\end{proof}

Now $\mathcal{M}(\theta_0,\theta_\infty)$ is, as algebraic variety, def\/ined as the quotient of
$\mathcal{A}(\theta_0,\theta_\infty) $ by $G$.
Since this is a~geometric quotient, the map $\mathcal{M}(\theta_0,\theta_\infty)(\mathbb{C})\rightarrow {\bf
S}(\alpha,\beta)$ is bijective.

\begin{Remarks}\label{R2.2}
1.~As in Section~\ref{Section1}, the Observation, one sees that, after scaling some of the variables, the two charts $ST_1$, $ST_0$ of
$\mathcal{M}(\theta_0,\theta_\infty)$, have the form $U_1\times T$, $U_2\times T$ with simply connected spaces
$U_1$, $U_2$.

{\it Define $M(\theta_0,\theta_\infty)=f^{-1}(1)$, where $f: \mathcal{M}(\theta_0,\theta_\infty)\rightarrow
T$ is the canonical morphism}.
Then $M(\theta_0,\theta_\infty)$ is simply connected.

Further $\tilde{\mathcal{M}}(\theta_0,\theta_\infty):= \mathcal{M}(\theta_0,\theta_\infty)\times_T\tilde{T}$ is the
universal covering of $\mathcal{M}(\theta_0,\theta_\infty)$.

2.~Let $q:\mathcal{M}(\theta_0,\theta_\infty)\rightarrow \mathbb{P}^1$ denote the morphism given by $q=-\frac{c_0}{c_1}$.

3.~In the cases $\alpha =\beta^{\pm 1}$, the space $\mathcal{A}(\theta_0,\theta_\infty)$ is def\/ined as before, and
including the assumption $(c_0,c_1)\neq (0,0)$.
The space $\mathcal{M}(\theta_0,\theta_\infty)$ is again the geometric and non singular quotient of
$\mathcal{A}(\theta_0,\theta_\infty)$ by $G$.
The canonical map $\mathcal{M}(\theta_0,\theta_\infty)(\mathbb{C})\rightarrow {\bf S}(\alpha,\beta)$ {\it is injective
and, in general, not surjective}.

Indeed, the open condition $(c_0,c_1)\neq (0,0)$ is valid for $(M,t)\in {\bf S}(\alpha,\beta)$ such that $M$ is
irreducible but may exclude certain reducible modules in ${\bf S}(\alpha,\beta)$.
We note that $c_0=c_1=0$ implies $\pm \frac{\theta_\infty}{2}=\frac{\theta_0}{2}-\epsilon$ with $\epsilon \in \{0,1\}$.
\end{Remarks}

\subsubsection[The Okamoto--Painlev\'e space $\tilde{\mathcal{M}}(\theta_0,\theta_\infty)$
for $\frac{\theta_0}{2}\pm \frac{\theta_\infty}{2}\not \in \mathbb{Z}$]{The Okamoto--Painlev\'e
space $\boldsymbol{\tilde{\mathcal{M}}(\theta_0,\theta_\infty)}$
for $\boldsymbol{\frac{\theta_0}{2}\pm \frac{\theta_\infty}{2}\not \in \mathbb{Z}}$}\label{section2.3.2}

The moduli space $\mathcal{M}(\theta_0,\theta_\infty)$ is replaced by $\tilde{\mathcal{M}}(\theta_0,\theta_\infty):= \mathcal{M}(\theta
_0,\theta_\infty)\times_T\tilde{T}$.
The bijections $\tilde{\bf S}(\alpha,\beta)\!\rightarrow \mathcal{R}(\alpha,\beta )\times \tilde{T}$ and
$\mathcal{M}(\theta_0,\theta_\infty)\!\rightarrow {\bf S}(\alpha,\beta)$ imply that the analytic morphism
$\tilde{\mathcal{M}}(\theta_0,\theta_\infty)\!\rightarrow \mathcal{R}(\alpha,\beta)\times \tilde{T}$ is bijective and
hence an analytic isomorphism.
As in Theorem~\ref{T1.3}, using arguments presented in~\cite{vdP2,vdP1,vdP-T} this implies the following result.

\begin{Theorem}
\label{T2.3}
Suppose that $\frac{\theta_0}{2}\pm \frac{\theta_\infty}{2}\not \in \mathbb{Z}$ $($equivalently $\alpha \neq \beta^{\pm1})$.
Then $\tilde{\mathcal{M}}(\theta_0,\theta_\infty)\rightarrow \tilde{T}$, provided with the foliation given by the
fibers of $\tilde{\mathcal{M}}(\theta_0,\theta_\infty)\rightarrow \mathcal{R}(\alpha,\beta)$ $($i.e., the isomonodromy
families$)$, is the Okamoto--Painlev\'e space corresponding to the equation ${\rm PIII(D_{6})}$, namely{\samepage
\begin{gather*}
q''=\frac{(q')^2}{q}-\frac{q'}{t}-\frac{4(\theta_0-1)}{t}+\frac{4\theta_\infty q^2}{t}+4q^3-\frac{4}{q}.
\end{gather*}
Moreover, this equation satisfies the Painlev\'e property.}
\end{Theorem}

\begin{observations}%\label{O2.4}
1.~The above formula dif\/fers slightly from the one given in~\cite[Section 4.5]{vdP-Sa}.
This is due to dif\/ferent choices of the standard matrix dif\/ferential operator.

2.~The transformation $t\mapsto -t$, $\theta_\infty \mapsto -\theta_\infty$ and $\theta_0\mapsto -\theta_0+2$ leaves the
family of matrix dif\/ferential operators invariant.
This has the consequence that a~solution $q(t)$ of ${\rm PIII(D_{6})}$ with parameters $\theta_0$ and $\theta_\infty$,
yields the solution $q(-t)$ of ${\rm PIII(D_{6})}$ with parameters $ -\theta_0+2$ and $-\theta_\infty$.
This can also be seen directly from the dif\/ferential equation.

3.~The solutions $q(t)$ of ${\rm PIII(D_6)}$ are in fact meromorphic functions in $\tilde{t}\in \tilde{T}=\mathbb{C}$.
Thus $Q(\tilde{t}):=q(e^{\tilde{t}})$ is well def\/ined and satisf\/ies the equation
\begin{gather*}
Q''=\frac{(Q')^2}{Q}-4(\theta_0-1)e^{\tilde{t}}+4\theta_\infty Q^2e^{\tilde{t}}+4Q^3e^{2\tilde{t}}-
\frac{4e^{2\tilde{t}}}{Q},
\qquad
\text{where}
\qquad
 '=\frac{d}{d\tilde{t}}.
\end{gather*}

4.~The space of initial conditions is analytically isomorphic to $\mathcal{R}(\alpha,\beta)$ and can also be identif\/ied
with $M(\theta_0,\theta_\infty)$.
Indeed, the extended Riemann--Hilbert isomorphism $\tilde{\mathcal{M}}(\theta_0,\theta_\infty) \rightarrow
\mathcal{R}(\alpha, \beta)\times \tilde{T}$ induces an analytic isomorphism $M(\theta_0,\theta_\infty)\rightarrow
\mathcal{R}(\alpha,\beta)$.
\end{observations}

\subsubsection{Verif\/ication of the formula in Theorem~\ref{T1.3}}
\label{Section2.3.3}

On the chart $ST_1$ of $\mathcal{M}(\theta_0,\theta_\infty)$ the matrix dif\/ferential operator has the form
$z\frac{d}{dz}+A=z\frac{d}{dz}+a_{-1}z^{-1}H+bE_1+(z-q)E_2$,
where $H=\begin{pmatrix}1&0 \\ 0&-1\end{pmatrix}$, $E_1=\begin{pmatrix}0&1\\ 0&0\end{pmatrix}$, $E_2=\begin{pmatrix}0&0\\1&0\end{pmatrix}$.
We will use $[H,E_1]=2E_2$, $[H,E_2]=-2E_2$, $[E_1,E_2]=H$. Further
\begin{gather*}
q:=-c_0,
b=b_{-2}z^{-2}+b_{-1}z^{-1}+b_0+b_1z,
\end{gather*}
where $b_{-2}\neq 0$ (by assumption),
\begin{gather*}
b_1=\frac{t^2}{4},
\qquad
b_0=q\frac{t^2}{4}+\frac{t\theta_\infty}{2},
\qquad
b_{-2}=a_{-1}+qb_{-1}+ t\left(\frac{\theta_0}{2}-\frac{1}{2}\right),
\\
a_{-1}^2-q\left(a_{-1}+qb_{-1}+t\left(\frac{\theta_0}{2}-\frac{1}{2}\right)\right)-\frac{t^2}{4}=0.
\end{gather*}
Now $q$, $a_{-1}$ are considered as (meromorphic) functions of $t$ and $A$ is a~matrix depending on $z$ and $t$.
The family $z\frac{d}{dz}+A$ is isomonodromic if and only if there is an operator of the form
$\frac{d}{dt}+B=\frac{d}{dt}+B_HH+B_1E_1+B_2E_2$ which commutes with $z\frac{d}{dz}+A$.
This is equivalent to $\frac{d}{dt}(A)=z\frac{d}{dz}(B)+[B,A]$.

Further $B_*=B_H$, $B_1$, $B_2$ are functions of $t$, $z$ and are supposed to have the form
$B_{*,-2}(t)z^{-2}+B_{*,-1}(t)z^{-1}+B_{*,0}(t)+B_{*,1}(t)z$.
The two operators commute if and only if $ a_{-1}'z^{-1}H+b'E_1-q'E_2$ is equal to
\begin{gather*}
\overset{\rm o}{B_H}H+\overset{\rm o}{B_1}E_1+\overset{\rm o}{B_2}E_2-\big[B_HH+B_1E_1+B_2E_2,a_{-1}z^{-1}H+bE_1+(z-q)E_2\big],
\end{gather*}
where $\overset{\rm o}{X}$ stands for $z\frac{d}{dz}(X)$ and $X':=\frac{d}{dt}(X)$.
On obtains the equations
\begin{alignat*}{3}
& (H)
\quad &&
a_{-1}'z^{-1}=\overset{\rm o}{B_H}-B_1(z-q)+B_2b,&
\\
& (E_1)
\quad&&
b'=\overset{\rm o}{B_1}-2B_Hb+2B_1a_{-1}z^{-1},&
\\
& (E_2)
\quad&&
-q'=\overset{\rm o}{B_2}+2B_H(z-q)-2B_2a_{-1}z^{-1}.&
\end{alignat*}

A Maple computation shows that this system of dif\/ferential equations for $q$, $a_{-1}$ is equivalent to
\begin{gather*}
q'=\frac{4a_{-1}-q}{t},\qquad a_{-1}'=\frac{4a_{-1}^2-t^2+q(t-a_{-1}-t\theta_0) +q^3t\theta_\infty +q^4t^2}{tq}.
\end{gather*}
The equation
\begin{gather*}
q''=\frac{(q')^2}{q}-\frac{q'}{t}-\frac{4(\theta_0-1)}{t}+\frac{4\theta_\infty
q^2}{t}+4q^3-\frac{4}{q}
\end{gather*}
follows by substitution.
Using the transformation $z\mapsto z^{-1}$ one f\/inds that $Q=\frac{1}{q}$ satisf\/ies the ${\rm PIII(D_6)}$ equation with
$\theta_0-1$ and $\theta_\infty$ interchanged:
\begin{gather*}
Q''=\frac{(Q')^2}{Q}-\frac{Q'}{t}-\frac{4\theta_\infty}{t}+\frac{4(\theta_0-1)Q^2}{t}+4Q^3-\frac{4}{Q}.
\end{gather*}

\subsubsection{After a~remark by Yousuke Ohyama}\label{Section2.3.4}

Let $(M,t)\in {\bf S}(\alpha,\beta)$ with $\alpha =e^{\pi i\theta_0}$, $\beta=e^{\pi i \theta_\infty}$ have the
property that $M$ is irreducible.
Consider the connection $(\mathcal{W},\nabla)$ with generic f\/iber $M$ and locally represented by
\begin{gather*}
z\frac{d}{dz}+\left(
\begin{matrix}
\frac{tz^{-1}+\theta_0}{2}&0\!\!\!\!
\\
0&-\frac{tz^{-1}+\theta_0}{2}
\end{matrix}
\right) \quad\! \text{at} \   z=0 \qquad \! \text{and} \! \qquad z\frac{d}{dz}+\left(
\begin{matrix}
\frac{tz+\theta_\infty}{2}&0\!\!\!\!
\\
0&-\frac{tz +\theta_\infty}{2}
\end{matrix}
\right) \quad\! \text{at} \ z=\infty.
\end{gather*}
The second exterior product of $(\mathcal{W},\nabla )$ is trivial and thus $\Lambda^2\mathcal{W}$ has degree 0.
Since $M$ is irreducible one has $\mathcal{W}\cong O(k)\oplus O(-k)$ with $k\in \{0,1\}$.

Suppose that the connection $\mathcal{W}$ has type $O(1)\oplus O(-1)$.
Then we can identify $\mathcal{W}$ with $O([0])B_1\oplus O(-[0])B_2$.
Put $D:=\nabla_{z\frac{d}{dz}}$.
Now $DB_1$ is not a~multiple of $B_1$ since $M$ is irreducible.
After multiplying $B_1$ with a~scalar, the matrix of $D$ with respect to the basis $B_1$, $B_2$
has the form $\begin{pmatrix}\alpha&\beta \\ z&-\alpha\end{pmatrix}$ with $\alpha =\alpha_{-1}z^{-1}+\alpha_0+\alpha_1z$,
$\beta =\beta_{-3}z^{-3}+\dots+\beta_1z$.
The base vector $B_2$ can be replaced by $B_2+hB_1$ with $h=h_0+h_{-1}z^{-1}+h_{-2}z^{-2}$.
The result is a~new representation of~$D$, namely
\begin{gather*}
\begin{pmatrix}1&h\\ 0&1\end{pmatrix}^{-1}\left\{z\frac{d}{dz}+\begin{pmatrix}\alpha&\beta \\ z&-\alpha\end{pmatrix}\right\}
\begin{pmatrix}1&h\\ 0&1\end{pmatrix}= z\frac{d}{dz}+\begin{pmatrix}\alpha -hz&2\alpha h-h^2z+z\frac{dh}{dz}+\beta\\z&-\alpha +hz\end{pmatrix}.
\end{gather*}
For unique $h_0$, $h_{-1}$ and at most two values of $h_{-2}$, the last operator is
\begin{gather*}
z\frac{d}{dz}+\begin{pmatrix}a_{-1}z^{-1}b\\ z-a_{-1}z^{-1}\end{pmatrix},
\qquad
\text{where}
\qquad
b=b_{-2}z^{-2}+b_{-1}z^{-1}+b_0+b_1z.
\end{gather*}
Let $e_1$, $e_2$ denote the new basis.
Then $\mathcal{V}=Oe_1\oplus O(-[0])e_2$ and the corresponding point $\xi\in \mathcal{M}(\theta_0,\theta_\infty)$
satisf\/ies $q(\xi)=0$.
The converse holds, too.
One observes that $q^{-1}(0)\subset \mathcal{M}(\theta_0,\theta_\infty)$ has two connected components, each one
isomorphic to $\mathbb{A}^1\times T$.
We note that the map $Q:=\frac{1}{q}$ can also be used in this context, since for a~monodromic family $Q$ satisf\/ies
a~${\rm PIII(D_6)}$ equation (see the end of Section~\ref{Section2.3.3}).

\looseness=-1
According to Malgrange, the locus where the bundle $\mathcal{W}$ is not free is the tau-divisor.
Thus we f\/ind that the tau-divisor coincides with the locus $Q^{-1}(\infty)\subset \mathcal{M}(\theta_0,\theta
_\infty)$.
The statement: `the tau-divisor coincides with $q^{-1}(\infty)$'
holds for PI, PII, PIII($\rm D_7$), PIII($\rm D_8$), PIV, too
(see~\cite{vdP2,vdP1,vdP-T}).

\subsection[The cases $\alpha=\beta^{\pm 1}$]{The cases $\boldsymbol{\alpha=\beta^{\pm 1}}$}

\subsubsection{Geometric quotients of the monodromy data}\label{Section2.4.1}

We use here the notation of Section~\ref{Section2.2}.

1.~If $\alpha =\beta\neq \pm 1$, then $\mathcal{R}(\alpha,\alpha)$ has a~singular point, namely
$(x_1,x_2,x_3)=(0,-\alpha,\alpha +\alpha^{-1})$.
The preimage in $\mathcal{T}(\alpha,\alpha )$ of this singular point consists of the tuples $(a_1,a_2,\ell_1,\dots,\ell
_4)$ such that the matrices $L$, $\operatorname{top}_0$ have the form $\begin{pmatrix}\ell_1 0 \\ \ell_3\ell_4\end{pmatrix}$,
$\begin{pmatrix}\alpha&0\\ \frac{a_1}{\alpha}&\frac{1}{\alpha}\end{pmatrix}$
or $\begin{pmatrix}\ell_1&\ell_2\\ 0&\ell_4\end{pmatrix}$, $\begin{pmatrix}\alpha&\alpha a_2 \\0&\frac{1}{\alpha}\end{pmatrix}$.
In particular, $\mathcal{R}(\alpha,\alpha)$ is {\it not a~geometric quotient}.

The remedy consists of replacing $\mathcal{T}(\alpha,\alpha)$ by $\mathcal{T}(\alpha,\alpha )^*$ which is the complement
of the closed subset of $\mathcal{T}(\alpha,\alpha)$ given by the equations $\ell_2=\ell_3=a_1=a_2=0$.
We claim that $\mathcal{T}(\alpha,\alpha )^*$ has a~nonsingular geometric quotient by the action of $\mathbb{G}_m\times
\mathbb{G}_m$.
A proof is obtained by writing $\mathcal{T}(\alpha,\alpha )^*$ as the union of the four af\/f\/ine open subsets $\ell_2\neq
0$, $\ell_3\neq 0$, $a_1\neq 0$ and $a_2\neq 0$.
On each of these subsets one explicitly computes the quotient by $\mathbb{G}_m\times \mathbb{G}_m$, which turns out to
be nonsingular and geometric.
Gluing these four quotients produces the required geometric quotient which will be denoted by $\mathcal{R}(\alpha,\alpha
)^*$.

Let ${\bf S}(\alpha,\alpha )^*$ the complement in ${\bf S}(\alpha,\alpha)$ of the set of the modules which are direct
sums and $\tilde{{\bf S}}(\alpha,\alpha )^*={\bf S}(\alpha,\alpha )^*\times_T\tilde{T}$.
Then the canonical map $\tilde{{\bf S}}(\alpha,\alpha)^* \rightarrow \mathcal{R}(\alpha,\alpha)^*\times \tilde{T}$ is
bijective.

Def\/ine the closed space $\mathcal{T}(\alpha,\alpha )^*_{\text{\rm red}}$ of $\mathcal{T}(\alpha,\alpha )^*$ by the condition that
the data is reducible.
This space has two irreducible components, given in terms of the matrices $L$, $\operatorname{top}_0$ by:

(a) $\begin{pmatrix}\ell_1&0\\ \ell_3&\ell_4\end{pmatrix}$,
$\begin{pmatrix}\alpha&0\\ \frac{a_1}{\alpha}&\frac{1}{\alpha}\end{pmatrix}$ with $\ell_1\ell_4=1$ and $(\ell_3,a_1)\neq 0$.
One easily verif\/ies that the map which sends $(L, \operatorname{top}_0)$ to $(\ell_3:a_1)\in \mathbb{P}^1$ is the geometric quotient.

(b) $\begin{pmatrix}\ell_1&\ell_2\\ 0&\ell_4\end{pmatrix}$, $\begin{pmatrix}\alpha&\alpha a_2\\ 0&\frac{1}{\alpha}\end{pmatrix}$
with $\ell_1\ell_4=1$ and $(\ell_2,a_2)\neq 0$.
One easily verif\/ies that the map which sends the $(L,\operatorname{top}_0)$ to $(\ell_2:a_2)\in \mathbb{P}^1$ is the geometric
quotient.

Therefore the `reducible locus' $\mathcal{R}(\alpha,\alpha )^*_{\text{\rm red}}$ (i.e., corresponding to reducible monodromy data)
is the union of two, not intersecting, projective lines.

2.~The case $\alpha =\beta^{-1}\neq \pm 1$ can be handled as in (1).
One f\/inds (with a~similar notation) a~geometric quotient $\mathcal{R}(\alpha,\alpha^{-1})^*$ of
$\mathcal{T}(\alpha,\alpha^{-1})^*$ and a~bijection $\tilde{{\bf S}}(\alpha,\alpha^{-1})^*\rightarrow
\mathcal{R}(\alpha,\alpha^{-1})^*\times \tilde{T}$.
Further $\mathcal{R}(\alpha,\alpha^{-1})^*_{\text{\rm red}}$ is the union of two, non intersecting, projective lines.

3.~$\alpha=\beta=1$.
The categorical quotient $\mathcal{R}(1,1)$ of $\mathcal{T}(1,1)$ has two singular points, namely
$(x_1,x_2,x_3)=(0,-1,2)$ and $(x_1,x_2,x_3)=(-1,0,2)$.
The preimage of the f\/irst singular point consists of the pairs $(L,\operatorname{top}_0)$
equal to $\left(\!\begin{pmatrix}\ell_1&\ell_2\\ 0&\ell_4\end{pmatrix},\begin{pmatrix}1&a_2\\0&1\end{pmatrix}\!\right)$
or to $\left(\!\begin{pmatrix}\ell_1&0\\ \ell_3&\ell_4\end{pmatrix},\begin{pmatrix}1&0\\ a_1&1\end{pmatrix}\!\right)$.
The~pre\-image of the second singular point consists of the pairs $(L,\operatorname{top}_0)$
equal to $\left(\!\begin{pmatrix}0&\ell_2\\ \ell_3&\ell_4\end{pmatrix},\begin{pmatrix}1&a_2\\ 0&1\end{pmatrix}\!\right)$
or to $\left(\!\begin{pmatrix}\ell_1&\ell_2\\ 0&\ell_4\end{pmatrix},\begin{pmatrix}1&0\\ a_1&1\end{pmatrix}\!\right)$.
Clearly, $\mathcal{R}(1,1)$ is not a~geometric quotient.

The locus of the points in $\mathcal{T}(1,1)$ which describe the monodromy data for modules in ${\bf S}(1,1)$ which are
direct sums is the union of the two closed sets $a_1=a_2=\ell_2=\ell_3=0$ and $a_1=a_2=\ell_1=\ell_4=0$.
Let $\mathcal{T}(1,1)^*\subset \mathcal{T}(1,1)$ denote the complement of this locus.
This set is the union of the six open subsets given by the inequalities $a_1\neq 0$, $a_2\neq 0$, $\ell_{12}\neq 0$,
$\ell_{13}\neq 0$, $\ell_{24}\neq 0$ and $\ell_{34}\neq 0$.
The group $\mathbb{G}_m\times\mathbb{G}_m$ acts on each of these open af\/f\/ine sets and the categorical quotient is
a~geometric quotient and is nonsingular.
Therefore the quotient $\mathcal{R}(1,1)^*$ of $\mathcal{T}(1,1)^*$, obtained by gluing the six quotients, is
a~geometric quotient and nonsingular.

{\it Example}.
The open subset $\ell_{12}\neq 0$ is def\/ined by the variables $\ell_1,\dots, \ell_4$, $a_1$, $a_2$
and relations $0=-\ell_{23}a_{12} +\ell_{24}a_1-\ell_{13}a_2$, $\ell_{14}-\ell_{13}-1=0$, $\ell_{12}\neq 0$.
Division by $\mathbb{G}_m\times \mathbb{G}_m$ is equivalent to the normalisation $\ell_1=\ell_2=1$.
Elimination of $\ell_4$ by $\ell_4=\ell_3+1$ yields the equation $\ell_3(a_1-a_2+a_1a_2)+a_1=0$.
This is a~nonsingular surface.
Using the projection $(\ell_3, a_1,a_2)\mapsto (a_1,a_2)$ one f\/inds that this surface is simply connected.

Similar computations lead to the statements:
{\it $\mathcal{R}(1,1)^*$ is a~nonsingular geometric quotient and is simply connected.
The natural map $\tilde{\bf S}(1,1)^{*}\rightarrow \mathcal{R}(1,1)^{*}\times \tilde{T}$ is a~well defined
bijection. The reducible locus $\mathcal{R}(1,1)^*_{\text{\rm red}}$ is the union of four, non intersecting, projective lines.}

4.~$\alpha =\beta =-1$.
The categorical quotient $\mathcal{R}(-1,-1)$ has two singular points, namely $(x_1,x_2,x_3)=(0,1,-2)$ and
$(x_1,x_2,x_3)=(-1,0,-2)$.
As in (3), one def\/ines $\mathcal{T}(-1,-1)^{*}$ and its geometric non singular quotient $\mathcal{R}(-1,-1)^{*}$.
The space $\mathcal{R}(-1,-1)^{*}$ contains four, non intersecting, projective lines.
These lines correspond to the reducible locus of $\mathcal{R}(-1,-1)^{*}$.
As in (3) one def\/ines $\tilde{{\bf S}}(-1,-1)^{*}$ and concludes:

{\it $\mathcal{R}(-1,-1)^*$ is a~nonsingular geometric quotient and is simply connected}.

{\it The natural map $\tilde{{\bf S}}(-1,-1)^{*}\rightarrow \mathcal{R}(-1,-1)^{*}\times \tilde{T}$ is a~bijection}.

\subsubsection[Reducible modules in ${\bf S}$]{Reducible modules in ${\bf S}$}

We use here the notation of Sections~\ref{Section2.2} and~\ref{Section2.4.1}.
\begin{observations}\label{O2.5}
Let $N\subset M$ be a~1-dimensional submodule, then $\mathbb{C}((z))\otimes N=\mathbb{C}((z))E_i$ and
$\mathbb{C}((z^{-1}))\otimes N=\mathbb{C}((z^{-1}))F_j$ with $i,j\in \{1,2\}$.
Since $N$ has no other singularities than $0$, $\infty$ one has $N=\mathbb{C}(z)n$ with $\delta (n)=(\frac{\pm tz^{-1}\pm
tz}{2}+d)n$, where $d\in \mathbb{C}$ is unique modulo $\mathbb{Z}$.
Indeed, $\delta (n)=\big(\frac{\pm tz^{-1}\pm tz}{2}+f\big)n$, where $f\in \mathbb{C}(z)$ has no poles at $0$ and $\infty$.
Using that $N$ has only singularities at $0$ and $\infty$, one can change the generator $n$ of $N$ such that $f$ is
a~constant $d$.
Any other base vector of $N$ with this property has the form $z^kn$ with $k\in \mathbb{Z}$.
Further $d\in \pm \theta_0/2+\mathbb{Z}$ and $d\in \pm \theta_\infty /2+\mathbb{Z}$ and hence $\alpha =\beta^{\pm
1}$.
\end{observations}

\begin{Proposition}[Reducible modules]\label{P2.6}

{\rm 1.}
A module $(M,t) \in {\bf S}$ is reducible if and only if there are $i,j\in \{1,2\}$
such that $a_i=0$, $b_j=0$ and
$L(\mathbb{C}e_i)=\mathbb{C}f_j$.

{\rm 2.}
Let $(M,t)\in {\bf S}$ be reducible, but not a~direct sum of two submodules of dimension one.
Then there are unique elements $\epsilon_1,\epsilon_2\in \{-1,1\}$,
a~complex number $d$, unique modulo $\mathbb{Z}$,
and a~polynomial $c=c_1z+c_0\neq 0$, unique up to multiplication by a~scalar, such that $M$ is represented by the matrix
differential operator
\begin{gather*}
z\frac{d}{dz}+\left(
\begin{matrix}
-\frac{\epsilon_1tz^{-1}+\epsilon_2tz}{2}-d&0
\\
c_1z+c_0&\frac{\epsilon_1tz^{-1}+\epsilon_2tz}{2}+d
\end{matrix}
\right).
\end{gather*}

{\rm 3.}
For $\alpha \neq \pm 1$, the reducible locus ${\bf S}(\alpha,\alpha )^*_{\text{\rm red}}$ of ${\bf S}(\alpha,\alpha )^*$ is
represented by the union of the two families in {\rm (2)} given by $e^{2\pi id}=\alpha$, $\epsilon_1=\epsilon_2=1$ and
$\epsilon_1=\epsilon_2=-1$.
Each of the two families is isomorphic to $\mathbb{P}^1\times T$, by sending the matrix differential operator to
$((c_1:c_0),t)\in \mathbb{P}^1\times T$.

The isomorphism $\tilde{{\bf S}}(\alpha,\alpha)^*_{\text{\rm red}}\rightarrow \mathcal{R}(\alpha,\alpha )^*_{\text{\rm red}}\times \tilde{T}$
yields two isomorphism $\mathbb{P}^1\times \tilde{T}\rightarrow \mathbb{P}^1\times \tilde{T}$.
These have the form $(p,\tilde{t})\mapsto (A(\tilde{t})p,\tilde{t})$ where $A(\tilde{t})$ is an automorphism of
$\mathbb{P}^1$, depending on~$\tilde{t}$.

{\rm 4.}
A similar result holds for $\alpha \neq \pm 1$ and $\alpha =\beta^{-1}$.

{\rm 5.}
${\bf S}(1,1)^{*}_{\text{\rm red}}$ is represented by the families $d\in \mathbb{Z}$ and the four possibilities of $\epsilon
_1,\epsilon_2$.
Then $\tilde{{\bf S}}(1,1)^{*}_{\text{\rm red}}$ identif\/ies with the disjoint union of four copies of $\mathbb{P}^1\times
\tilde{T}$.
The same holds for $\mathcal{R}(1,1)^{*}_{\text{\rm red}}\times \tilde{T}$.
The isomorphism $\tilde{{\bf S}}(1,1)^{*}_{\text{\rm red}}\rightarrow \mathcal{R}(1,1)^{*}_{\text{\rm red}}\times \tilde{T}$ yields four
isomorphism $\mathbb{P}^1\times \tilde{T}\rightarrow \mathbb{P}^1\times \tilde{T}$.
These have the form described in {\rm (3)}.

{\rm 6.}
The case $\alpha =\beta =-1$ is similar to case {\rm (5)}.
\end{Proposition}

\begin{proof}
1.~This follows from Observations~\ref{O2.5} and the statement that a~dif\/ferential module over $\mathbb{C}(z)$ is
determined by its monodromy data (i.e., ordinary monodromy, Stokes matrices and links) and the formal classif\/ication of
the singular points (see~\cite[Theorem~1.7]{vdP-Sa}).

2.~For convenience we consider the case $i=j=1$ of (1).
Using the above Observation, one f\/inds that $M$ has a~basis $m_1$, $m_2$ such that $z\partial (m_2)=am_2$ and $z\partial
(m_1)=-am_1+fm_2$ with $a:=\frac{tz^{-1}+tz}{2}+d$ and $f\in \mathbb{C}(z)$.
If we f\/ix $d$, then $m_2$ is unique up to multiplication by a~scalar.
Further, $m_1$ is unique up to a~transformation $m_1\mapsto \lambda m_1+hm_2$ with $\lambda \in \mathbb{C}^*$, $h\in
\mathbb{C}(z)$.
This transformation changes $f$ into $\lambda f+2ah +zh'$.

We start considering the subgroup of transformations with $\lambda =1$ and $h\in \mathbb{C}(z)$.
For a~suitab\-le~$h$ the term $f_1:=f+2ah+zh'$ has in $\mathbb{C}^*$ at most poles of order one.
A pole of order one of~$f_1$ in~$\mathbb{C}^*$ cannot disappear by a~transformation of the form under consideration.
Since~$M$ has only singularities at $0$ and $\infty$ we conclude that $f_1\in \mathbb{C}[z,z^{-1}]$.
For suitable $h\in \mathbb{C}[z,z^{-1}]$ the term $c:=f_1+2ah+zh'$ is a~polynomial of degree ${\leq}1$ and is ${\neq}0$,
by assumption.
For any $h\in \mathbb{C}(z)$, $h\neq 0$ the term $c+2ah+zh'$ is not a~polynomial of degree ${\leq}1$.
This yields a~unique $c$ for this subgroup of transformations.
Finally, the transformation $m_1\mapsto \lambda m_1$ shows that~$c$ is unique up to multiplication by a~scalar.

Cases (3)--(6) are consequences of the above computation of the $\mathcal{R}(\alpha,\beta )^*_{\text{\rm red}}$.
\end{proof}

\subsubsection[The reducible connections in $\mathcal{M}(\theta_0,\theta_\infty)$]{The reducible
connections in $\boldsymbol{\mathcal{M}(\theta_0,\theta_\infty)}$}
\label{Section2.4.3}

The image of the injective map $\mathcal{M}(\theta_0,\theta_\infty)\rightarrow {\bf S}(\alpha,\beta)^*$, where
$\alpha=e^{\pi i \theta_0}$, $\beta=e^{\pi i \theta_\infty}$, will be denoted by ${\bf S}(\theta_0,\theta_\infty)$.
By Remarks~\ref{R2.2} part (3), this image contains the irreducible modules $(M,t)$.
For $\alpha\neq \beta^{\pm 1}$, ${\bf S}(\theta_0, \theta_\infty)$ is equal to ${\bf S}(\alpha,\beta)={\bf
S}(\alpha,\beta)^*$.
Now we consider the other cases.

\begin{Proposition}%\label{P2.7}
Suppose $\alpha=\beta \neq \pm 1$.
Then ${\bf S}(\theta_0,\theta_\infty)$ consists of:
\begin{enumerate}\itemsep=0pt
\item[{\rm (a)}] the irreducible elements $(M,t)$;

\item[{\rm (b)}] the $(\epsilon_1,\epsilon_2)=(1,1)$ family of reducible modules $\Leftrightarrow$ $\theta_0- \theta_\infty
\geq 2$;

\item[{\rm (c)}] the $(\epsilon_1,\epsilon_2)=(-1,-1)$ family of reducible modules $\Leftrightarrow$ $\theta_0- \theta_\infty
\leq 0$.
\end{enumerate}
\end{Proposition}
\begin{proof}
Part (a) is known.
Consider an element $(M,t)$ in the $(\epsilon_1,\epsilon_2)=(1,1)$ family of reducible modules.
Let the connection $(\mathcal{V},\nabla)$, corresponding to $(M,t)$, be def\/ined as in Section~\ref{Section2.3.1}.
Then $\mathcal{V}$ can be identif\/ied with the vector bundle $O(k[0])e_1+O((-k-1)[0])e_2$ for some integer $k\geq 0$.
Then $(M,t)$ lies in ${\bf S}(\theta_0,\theta_\infty)$ $\Leftrightarrow$ $k=0$.

Consider the case $k>0$.
Then $\nabla_{z\frac{d}{dz}}e_1=(a_{-1}z^{-1}+a_0+a_1z)e_1$ and
$\nabla_{z\frac{d}{dz}}(z^{-k}e_1)=(a_{-1}z^{-1}+a_0+a_1z-k)(z^{-k}e_1)$.

Comparing with the prescribed local operator at $z=0$ yields the possibilities: $a_{-1}z^{-1}+a_0$ equals (i)
$\frac{t}{2}z^{-1}+\frac{\theta_0}{2}$ or (ii) $-\frac{t}{2}z^{-1}-\frac{\theta_0}{2}+1$.

Comparing with the prescribed local operator at $z=\infty$ yields the possibilities: $a_0+a_1z$ equals (A)
$\frac{t}{2}z+\frac{\theta_\infty}{2}$ or (B) $-\frac{t}{2}z-\frac{\theta_\infty}{2}$.
Combining one obtains
\begin{alignat*}{7}
& (i),\quad && (A)
\qquad &&
a_{-1}=\frac{t}{2},\qquad && a_1=\frac{t}{2},\qquad && \theta_0=\theta_\infty -2k,\qquad && (\epsilon_1,\epsilon_2)=(1,1),&
\\
& (i),\quad && (B)
\qquad &&
a_{-1}=\frac{t}{2},\qquad && a_1=-\frac{t}{2},\qquad && \theta_0=-\theta_\infty -2k,\qquad && (\epsilon_1,\epsilon_2)=(1,-1),&
\\
& (ii),\quad && (A)
\qquad &&
a_{-1}=-\frac{t}{2},\qquad &&a_1=\frac{t}{2},\qquad && \theta_0=-\theta_\infty +2k +2,\qquad && (\epsilon_1,\epsilon_2)=(-1,1),&
\\
& (ii),\quad && (B)
\qquad &&
a_{-1}=-\frac{t}{2},\qquad && a_1=-\frac{t}{2},\qquad && \theta_0=\theta_\infty +2k+2,\qquad && (\epsilon_1,\epsilon_2)=(-1,-1).&
\end{alignat*}
Since $\alpha =\beta \neq \pm 1$, reducible modules of types $(1,-1)$ and $(-1,1)$ are not present in ${\bf
S}(\alpha,\alpha )^*$.
Now we consider the presence of reducible modules of type $(1,1)$ in ${\bf S}(\theta_0,\theta_\infty)$.The condition
$\theta_0- \theta_\infty\geq 0$ is necessary because of $(i)$,~$A$.

Consider the case $\theta_0=\theta_\infty=2d$.
A reducible module of type $(1,1)$ yields a~connection on $\mathcal{W}=Of_1\oplus Of_2$ with the local data
$z\frac{d}{dz}+\begin{pmatrix}-(\frac{t}{2}z^{-1}+d)
&0\\ *&\frac{t}{2}z^{-1}+d\end{pmatrix}$ at $z=0$ and $z\frac{d}{dz}+\begin{pmatrix}-(\frac{t}{2}z+d)
&0\\ *
&\frac{t}{2}z+d\end{pmatrix}$ at $z=\infty$.
Now $\mathcal{V}=O(-[0])f_1\oplus Of_2=Oe_1\oplus O(-[0])e_2$, with $e_1=f_2$, $e_2=f_1$, has the required local data and
the matrix of $\nabla_{z\frac{d}{dz}}$ with respect to the basis $e_1$, $e_2$ is $\begin{pmatrix}\omega&*\\ 0&-\omega\end{pmatrix}$.
Thus $c_0=c_1=0$ and the reducible modules of type $(1,1)$ are not present according to the construction of
$\mathcal{M}(\theta_0,\theta_\infty)$.

Consider the case $\theta_0=\theta_\infty +2$.
The standard form of Proposition~\ref{P2.6}.2 for type $(1,1)$ belongs to $\mathcal{M}(2d+2,2d)$.
Further, see Remarks~\ref{R2.2}.3, $(c_0,c_1)\neq (0,0)$ holds for $\theta_0-\theta_\infty >2$.
This proves~(b).
The proof of~(c) is similar.
\end{proof}

Similarly one shows:
$\alpha =\beta^{-1}\neq \pm 1$.
Then ${\bf S}(\theta_0,\theta_\infty) $ consists of:
\begin{enumerate}\itemsep=0pt
\item[{\rm (a)}] the irreducible elements $(M,t)$;

\item[{\rm (b)}] the $(\epsilon_1,\epsilon_2)=(1,-1)$ family of reducible modules $\Leftrightarrow$ $\theta_0 +\theta_\infty
\geq 2$;

\item[{\rm (c)}] the $(\epsilon_1,\epsilon_2)=(-1,1)$ family of reducible modules $\Leftrightarrow$ $\theta_0 + \theta_\infty
\leq 0$.
\end{enumerate}

$\alpha =\beta=\pm 1$.
Then ${\bf S}(\theta_0,\theta_\infty) $ consists of:
\begin{enumerate}\itemsep=0pt
\item[{\rm (a)}] the irreducible elements $(M,t)$;

\item[{\rm (b)}] the $(\epsilon_1,\epsilon_2)=(1,1)$ family of reducible modules $\Leftrightarrow$ $\theta_0 - \theta_\infty\geq 2$;

\item[{\rm (c)}] the $(\epsilon_1,\epsilon_2)=(1,-1)$ family of reducible modules $\Leftrightarrow$ $\theta_0 +\theta_\infty\geq 2$;

\item[{\rm (d)}] the $(\epsilon_1,\epsilon_2)=(-1,1)$ family of reducible modules $\Leftrightarrow$ $\theta_0 +\theta_\infty\leq 0$;

\item[{\rm (e)}] the $(\epsilon_1,\epsilon_2)=(-1,-1)$ family of reducible modules $\Leftrightarrow$ $\theta_0 -\theta_\infty\leq 0$.
\end{enumerate}

\subsubsection[$\mathcal{M}(\theta_0, \theta_\infty)$ for $\alpha =\beta^{\pm 1}$]{$\boldsymbol{\mathcal{M}(\theta_0, \theta_\infty)}$
for $\boldsymbol{\alpha =\beta^{\pm 1}}$}\label{section2.4.4}

Let $\mathcal{R}(\theta_0.\theta_\infty)$ denote the open
subspace of $\mathcal{R}(\alpha, \beta)^*$ which corresponds to the subset ${\bf S}(\theta_0,\theta_\infty)$ of ${\bf
S}^*(\alpha,\beta)$, def\/ined in Section~\ref{Section2.4.3}.
By Sections~\ref{Section2.4.1}--\ref{Section2.4.3}, the extended Riemann--Hilbert
morphism $\tilde{\mathcal{M}}(\theta_0,\theta_\infty)\rightarrow\mathcal{R}(\theta_0,\theta_\infty)\times \tilde{T}$
is a~well def\/ined analytic isomorphism.
This has as consequence:

{\it Theorem~{\rm \ref{T2.3}} holds for the cases $\frac{\theta_0}{2}\pm \frac{\theta_\infty}{2}\in \mathbb{Z}$ with
$\mathcal{R}(\alpha,\beta)$ replaced by $\mathcal{R}(\theta_0,\theta_\infty)$}.

\subsubsection{Isomonodromy for reducible connections}

The f\/ibers of the locally def\/ined map ${\bf S}(\alpha,\alpha^{\pm 1})^*_{\text{\rm red}}\rightarrow \mathcal{R}(\alpha,\alpha^{\pm
1})^*_{\text{\rm red}}$ are the isomonodromy families of reducible modules.
As a~start, we consider the reducible familly of type $(\epsilon_1,\epsilon_2)=(1,1)$ lying in $\mathcal{M}(2d+2,2d)$.
This family is represented by $ z\frac{d}{dz}+\left(
\begin{matrix}
-\frac{tz^{-1}+tz}{2}-d&0
\\
z-q &\frac{tz^{-1}+tz}{2}+d
\end{matrix}
\right).
$ For an isomonodromy subfamily of this, $q$ is a~function of $t$ and the Stokes data at $0$ and $\infty$ and the link
are f\/ixed.
Isomonodromy is equivalent to the statement that the above matrix dif\/ferential operator commutes with an operator of the
form $\frac{d}{dt}+B_{-1}z+B_0+B_1z$, where the tracefree $2\times 2$ matrices $B_{-1}$, $B_0$, $B_1$ depend on $t$ only.
This leads to the equation
\begin{gather*}
\left(
\begin{matrix}
-\frac{z^{-1}+z}{2}&0
\\
-q'&\frac{z^{-1}+z}{2}
\end{matrix}
\right)=
-B_1z^{-1}+B_1z\\
\hphantom{\left(
\begin{matrix}
-\frac{z^{-1}+z}{2}&0
\\
-q'&\frac{z^{-1}+z}{2}
\end{matrix}
\right)=}{}
-\left[B_{-1}z^{-1}+B_0+B_1z, \left(
\begin{matrix}
-\frac{tz^{-1}+tz}{2}-d&0
\\
z-q&\frac{z^{-1}+z}{2}+d
\end{matrix}
\right)\right].
\end{gather*}
A computation yields the equation $q'=-2q^2-\frac{4d-1}{t}q-2$.
The solutions of this equation have the form $\frac{1}{2}\frac{y'}{y}$, where $y$ is a~non zero solution of the Bessel
equation $y''+\frac{4d-1}{t}y'+4y=0$.
One obtains in a~similar way for an isomonodromic family of reducible modules of type $(\epsilon_1,\epsilon_2)$ the equation
\begin{gather*}
q'=-2\epsilon_2q^2-\frac{4d-1}{t}q-2\epsilon_1.
\end{gather*}
The solutions are $q=\frac{\epsilon_2}{2}\frac{y'}{y}$ where $y$ is a~solution of the Bessel equation
$y''+\frac{4d+1}{t}y'+4\epsilon_1\epsilon_2y=0$.
These equations are consistent with the formula of Theorem~\ref{T2.3}
\begin{gather*}
q''=\frac{(q')^2}{q}-\frac{q'}{t}-\frac{4(\theta_0-1)}{t}+\frac{4\theta_\infty q^2}{t}+4q^3-\frac{4}{q}
\end{gather*}
for isomonodromic families in $\mathcal{M}(\theta_0,\theta_\infty)$.
According to~\cite{Oh} we found in this way all Riccati solutions for ${\rm PIII(D_6)}$, up to the action of the
B\"acklund transformations.

\begin{Remark}%\label{R2.8}
The assumption that a~function $q$ satisf\/ies two distinct ${\rm PIII(D_6)}$ equations leads to $q^4=1$.
Thus we found the algebraic solutions $q=\pm 1$ for $\theta_\infty =\theta_0-1$ and $q=\pm i$ for $-\theta_\infty
=\theta_0-1$.
According to~\cite{Oh} these are all the algebraic solutions of ${\rm PIII(D_6)}$, up to the action of the B\" acklund
transformations.
\end{Remark}

\section[B\"acklund transformations for ${\rm PIII(D_6)}$]{B\"acklund transformations for $\boldsymbol{\rm PIII(D_6)}$}

\subsection[Automorphisms of $\bf S$]{Automorphisms of ${\bf S}$}

We start with a~table of generators for the group ${\rm Aut}({\bf S})$ of `natural' automorphism of $\bf S$, in terms of
their action on the parameters $\alpha$, $\beta$ and on $t$, $z$.
\begin{center}
\begin{tabular}{|c|c|c|c|c|}
\hline
&$\alpha$ &$\beta$ &$t$&$z$
\\
\hline
$\sigma_1$&$\alpha^{-1}$&$\beta^{-1}$&$-t$&$z$\tsep{1pt}
\\
\hline
$\sigma_2$&$-\alpha$&$-\beta $&$t$&$z$
\\
\hline
$\sigma_3$&$\alpha$&$\beta^{-1}$&$it$&$iz$\tsep{1pt}
\\
\hline
$\sigma_4$&$\beta $&$\alpha$&$t$&$z^{-1}$\tsep{1pt}
\\
\hline
\end{tabular}
\end{center}

These generators are def\/ined as follows.

1.~$\sigma_1:(M,t)\mapsto (M,-t)$.
This induces bijections ${\bf S}(\alpha,\beta )\rightarrow {\bf S}(\alpha^{-1},\beta^{-1})$.
Indeed, the basis vectors $e_1$, $e_2$ of $V(0)$ are interchanged and the same holds for the basis $f_1$, $f_2$ of
$V(\infty)$.

2.~Def\/ine the dif\/ferential module $N=\mathbb{C}(z)b$ by $\delta b=\frac{1}{2}b$.
Then $\sigma_2: (M,t)\mapsto (M\otimes N,t)$.
Since $\Lambda^2(M\otimes N)=N^{\otimes 2}$ is the trivial module, $(M\otimes N,t)$ belongs to $\bf S$.
Let $E_1$, $E_2$ be a~basis of $\mathbb{C}((z))\otimes M$ such that $\delta E_1=-\frac{tz^{-1}+\theta_0}{2}E_1$ and
$\delta E_2=\frac{tz^{-1}+\theta_0}{2}E_2$.
Then the formal module $\mathbb{C}((z))\otimes (M\otimes N)$ has basis $E_1\otimes b$, $E_2\otimes b$ and $\delta
(E_1\otimes b)=(-\frac{tz^{-1}+\theta_0}{2}+\frac{1}{2})(E_1\otimes b)$ and similarly $\delta (E_2\otimes
b)=(\frac{tz^{-1}+\theta_0}{2}+\frac{1}{2})(E_2\otimes b)$.
Thus $e^{\pi i (\theta_0 +1)}=-\alpha$ is the eigenvalue of the formal monodromy at $z=0$.
The same argument shows that $-\beta$ is the eigenvalue of the formal monodromy at $z=\infty$.

3.~Let $\phi$ be a~$\mathbb{C}$-linear automorphism of the f\/ield $\mathbb{C}(z)$, such that $z\frac{d}{dz}\circ
\phi=\mu\cdot \phi \circ z\frac{d}{dz}$ for some $\mu \in \mathbb{C}^*$.
There are two possibilities:

(a)~$\phi (z)=cz$ with $c\in \mathbb{C}^*$, $\mu =1$,
(b)~$\phi (z)=cz^{-1}$ with $c\in \mathbb{C}^*$, $\mu =-1$.

For $(M,t)\in {\bf S}$ one considers $(\mathbb{C}(z)\otimes_{\phi}M,\dots)$.
As additive group $\mathbb{C}(z)\otimes_{\phi}M$ is identif\/ied with $M$.
Its structure as vector space is given by the new scalar multiplication $f*m:=\phi (f)m$.
In case (a), the dif\/ferential structure is given by the original $\delta$ and in case~(b) the dif\/ferential structure is
given by $-\delta$.

(a)~$\phi (z)=cz$.
The formal local module $\mathbb{C}((z))\otimes M$ with basis $E_1$, $E_2$ and $\delta E_1=-\frac{tz^{-1}+\theta
_0}{2}E_1$, $\delta E_2=\frac{tz^{-1}+\theta_0}{2}E_2$ is transformed into $\mathbb{C}((z))\otimes_{\phi}M$.
Now $\delta E_1= -\frac{tcz^{-1}+\theta_0}{2}*E_1$ and $ \delta E_2=\frac{tcz^{-1}+\theta_0}{2}*E_2$.
The basis $F_1$, $F_2$ of $\mathbb{C}((z^{-1}))\otimes M$ with $\delta F_1=-\frac{tz+\theta_\infty}{2}F_1$ and $\delta
F_2=\frac{tz+\theta_\infty}{2}F_2$, yields for the new structure the formulas $\delta F_1=-\frac{tc^{-1}z+\theta
_\infty}{2}*F_1$ and $\delta F_2=\frac{tc^{-1}z+\theta_\infty}{2}*F_2$.

The condition $tc=\pm tc^{-1}$ implies that $c^4=1$.
We def\/ine $\sigma_3$ by $\phi (z)=iz$.
This yields the new module $(\mathbb{C}(z)\otimes_\phi M,it)$
with new $\alpha$ equal to $e^{\pi i \theta_0}$ and new
$\beta$ equal to $e^{ -\pi i \theta_\infty}$.

(b)~We def\/ine $\sigma_4$ by $\phi (z)=z^{-1}$.
The new module $(\mathbb{C}(z)\otimes_\phi M,t)$ has new $\alpha =e^{\pi i \theta_\infty}$ and new $\beta =e^{\pi i
\theta_0}$.

{\it Comments.}
Using the f\/irst three columns of the table we will consider ${\rm Aut}({\bf S})$ as an automorphism group of the
algebraic variety $\mathbb{C}^*\times \mathbb{C}^*\times \mathbb{C}^*$.
A straightforward computation shows that ${\rm Aut}({\bf S})$ is the product $\langle \sigma_1\rangle \times \langle \sigma_2\rangle \times
\langle \sigma_3,\sigma_4\rangle $, where $\langle \sigma_1\rangle $ and $\langle \sigma_2\rangle $ have order two and $\langle \sigma_3,\sigma_4\rangle $ is the dihedral
group $D_4$ of order eight.

{\it The action of ${\rm Aut}({\bf S})$ on the monodromy data}.

1.~$\sigma_1$.
The map $t\mapsto -t$ has as consequence that the basis vectors $e_1$, $e_2$ of $V(0)$ are permuted and the same holds for
$V(\infty)$.
The singular directions and the Stokes maps do not change.
The new topological monodromies are $\left(
\begin{matrix}
\frac{1+a_1a_2}{\alpha}&\frac{a_1}{\alpha}
\\
\alpha a_2&\alpha
\end{matrix}
\right)$ at $z=0$ and $\left(
\begin{matrix}
\frac{1+b_1b_2}{\beta}&\frac{b_1}{\beta}
\\
\beta b_2&\beta
\end{matrix}
\right)$ at $z=\infty$.
The matrix of the link $L$ is now $\begin{pmatrix}\ell_4&\ell_3\\ \ell_2&\ell_1\end{pmatrix}$.
The matrix relation remains the same.
This amounts to the change $\alpha \mapsto \alpha^{-1}$, $\beta \mapsto \beta^{-1}$, $a_1\leftrightarrow a_2$,
$b_1\leftrightarrow b_2$, $\ell_1\leftrightarrow \ell_4$, $\ell_2\leftrightarrow \ell_3$.
This induces an automorphism of~$\mathcal{R}$ given by the formula $(\alpha, \beta, x_1,x_2,x_3)\mapsto (\alpha
^{-1},\beta^{-1},x_1\alpha^{-1}\beta^{-1},x_2\alpha^{-1}\beta^{-1},x_3)$.

2.~$\sigma_2$.
The map $(M,t) \mapsto (M\otimes N,t)$ has as consequence that the formal monodromies at $z=0$ and $z=\infty$ are
multiplied by $\begin{pmatrix}-1&0\\ 0&-1\end{pmatrix}$ (and thus $\alpha \mapsto -\alpha$, $\beta \mapsto -\beta$).
The Stokes matrices do not change.
The same holds for the link $L$.
The induced automorphism of~$\mathcal{R}$ is given by the formula $(\alpha,\beta,x_1,x_2,x_3)\mapsto (-\alpha,
-\beta,x_1,-x_2,-x_3)$.

3.~$\sigma_3$.
The ef\/fect of this transformation is: the singular directions change over $\frac{\pi}{2}$; the basis vectors of $V(0)$
are permuted; the basis of $V(\infty)$ is unchanged; the Stokes map and the link remain the same, however the
corresponding matrices change.
The induced automorphism of~$\mathcal{R}$ is given by $(\alpha,\beta,x_1,x_2,x_3)\mapsto (\alpha^{-1},\beta,\alpha
^{-1}x_1,\alpha^{-1}x_2,x_3)$.

4.~$\sigma_4$.
The spaces $V(0)$ and $V(\infty)$ are interchanged; the other data are unchanged.
The induced automorphism of~$\mathcal{R}$ is given by $(\alpha,\beta,x_1,x_2,x_3)\mapsto (\beta,\alpha,x_1,x_2,x_3)$.

\subsection{B\"acklund transformations}

Def\/ine the map $exp:\mathbb{C}^3\rightarrow (\mathbb{C}^*)^3$ by $(\theta_0,\theta_\infty,\tilde{t})\mapsto (\alpha,
\beta,t):=(e^{\pi i \theta_0},e^{\pi i \theta_\infty}, e^{\tilde{t}})$.
We will def\/ine the group $B({\bf S})$ of the {\it B\"acklund transformations} as the af\/f\/ine automorphisms of
$\mathbb{C}^3$ which respect the equivalence relation def\/ined by $exp$ and which map to elements of ${\rm Aut}({\bf
S})$.
By def\/inition there is an exact sequence of groups $0\rightarrow B({\bf S})_0\rightarrow B({\bf S})\rightarrow {\rm
Aut}({\bf S})\rightarrow 1$, where $B({\bf S})_0$ is the group of af\/f\/ine transformations of $\mathbb{C}^3$, generated
by: $B_1:(\theta_0,\theta_\infty,\tilde{t})\mapsto (2+\theta_0, \theta_\infty, \tilde{t})$, $B_2:(\theta_0,\theta
_\infty,\tilde{t})\mapsto (\theta_0,2+ \theta_\infty, \tilde{t})$ and $B_3:(\theta_0,\theta_\infty,\tilde{t})\mapsto
(\theta_0, \theta_\infty,2\pi i + \tilde{t})$.

The aim is to give each B\"acklund transformation the interpretation of a~morphism between the various moduli spaces
$\tilde{\mathcal{M}}(*,*)$, preserving the foliations by isomonodromic families, and to compute the ef\/fect on solutions
of ${\rm PIII(D_6)}$.
First we investigate the group $B({\bf S})$.

The af\/f\/ine map $B_3$ is not considered in the literature.
Its action on $\tilde{\mathcal{M}}(\theta_0,\theta_\infty)$ is obvious since $\tilde{\mathcal{M}}(\theta_0,\theta
_\infty)=\mathcal{M}(\theta_0,\theta_\infty)\times_T\tilde{T}$ and $\tilde{T}=\mathbb{C}\rightarrow T=\mathbb{C}^*$
is the map $\tilde{t}\mapsto e^{\tilde{t}}$.
The ef\/fect of $B_3$ on solutions of ${\rm PIII(D_6)}$ is far from obvious.
A solution is a~function $q(\tilde{t})=\text{``}q(e^{\tilde{t}})\text{''}$.
Since the equation depends only on $t$, the function $q(2\pi i +\tilde{t})$ satisf\/ies the same equation.
It seems that no formula for $q(2\pi i +\tilde{t})$ in terms of $q(\tilde{t})$ and its derivative is present in the
literature.

Generators for $B({\bf S})$ are given by their action on $\mathbb{C}^3$ and the variables $t$, $z$.
\begin{center}
\begin{tabular}{|c|c|c|c|c|c|}
\hline
&$\theta_0$&$\theta_\infty$&$\tilde{t}$&$ t$ &$z$\tsep{1pt}
\\
\hline
$s_1$ &$2-\theta_0 $&$-\theta_\infty $&$ \pi i +\tilde{t}$ &$-t $&$ z $\tsep{1pt}
\\
\hline
$s_2$ &$ 1+\theta_0$&$1+\theta_\infty $&$\tilde{t} $ &$ t$ &$ z$\tsep{1pt}
\\
\hline
$s_3$ &$\theta_0 $&$-\theta_\infty $&$ \frac{\pi i}{2}+\tilde{t}$ &$it $&$iz $\tsep{1pt}
\\
\hline
$s_4$ &$\theta_\infty $&$\theta_0 $&$\tilde{t} $ &$t $&$z^{-1} $\tsep{1pt}
\\
\hline
\end{tabular}
\end{center}
These elements generate $B({\bf S})$ because: $s_*$ is mapped to $\sigma_*$ for $*=1,2,3,4$;
$B_3=s_1^2=s_3^4$; $B_1=s_3s_1^{-1}s_4s_3s_4$ and $B_2=B_1^{-1}s_2^2$.

The group $\langle B_3\rangle $, generated by $B_3$, is isomorphic to $\bf Z$ and lies in the center of $B({\bf S})$.
Put $\overline{B({\bf S})}=B({\bf S})/\langle B_3\rangle $ and let $\overline{s}_*$ denote the image of $s_*$ in this quotient.
Then $\overline{s}^2_3$ has order two and lies in the center of $\overline{B({\bf S})}$.

We compare this with Okamoto's paper~\cite{O4}.
The group of the B\"acklund transformations $B$ of equation ${\rm PIII'(D_6)}$
\begin{gather*}
\frac{d^2Q}{dx^2}=\frac{1}{Q}\left( \frac{dQ}{dx} \right)^2-\frac{1}{x}\frac{dQ}{dx} +\frac{Q^2(\gamma Q+\alpha)}{4x^2}
+\frac{\beta}{4x}+\frac{\delta}{4Q}
\end{gather*}
is computed to be the af\/f\/ine Weyl group of type $B_2$.
The substitution $x=t^2$, $Q=tq$ transforms the equation into
\begin{gather*}
\frac{d^2q}{dt^2}=\frac{1}{q}\left(\frac{dq}{dt}\right)^2-\frac{1}{t}\frac{dq}{dt}+\frac{\alpha q^2+\beta}{t}+\gamma
q^3+\frac{\delta}{q},
\end{gather*}
and thus in our notation $\alpha =4\theta_\infty$, $\beta =-4(\theta_0 -1)$, $\gamma =4$, $\delta =-4$.
In a~sense, ${\rm PIII(D_6)}$ is a~degree two covering of ${\rm PIII'(D_6)}$.
It can be seen that there is a~surjective homomorphism $\overline{B({\bf S})}\rightarrow B$ with kernel $\langle \overline{s}_3^2\rangle $.

Our approach using moduli spaces explains the B\"{a}cklund transformations presented in~\cite{MCB}.
In contrast to this, the new transformations of~\cite{W} do not seem to have a~simple modular interpretation.

\subsection{B\"acklund transformations of the moduli spaces}

Let $s\in B({\bf S})$ have image $\sigma \in {\rm Aut}({\bf S})$.
Choose $\alpha$, $\beta$ and write $\alpha '=\sigma (\alpha)$, $\beta '=\sigma (\beta)$.
Choose $\theta_0$, $\theta_\infty$, $\theta_0'$, $\theta_\infty'$ such that $e^{\pi i \theta_0}=\alpha,\dots,
e^{\pi i\theta_\infty '}=\beta '$ and $s(\theta_0)=\theta_0'$, $s(\theta_\infty)=\theta_\infty '$.
Now $\sigma$ induces a~bijection ${\bf S}(\alpha,\beta)\stackrel{\sigma}{\rightarrow}{\bf S}(\alpha',\beta ')$ and consider
\begin{gather*}
\mathcal{M}(\theta_0,\theta_\infty)\rightarrow {\bf S}(\alpha,\beta)\stackrel{\sigma}{\rightarrow} {\bf
S}(\alpha',\beta ')\leftarrow \mathcal{M}(\theta_0',\theta_\infty ').
\end{gather*}
The f\/irst and the last arrow are injective.
Their images contain the locus of the irreducible mo\-du\-les and the loci of the $(\epsilon_1,\epsilon_2)$-reducible
modules depending on the $\theta_0,\dots,\theta_\infty'$ (see Section~\ref{Section2.4.3}).
Thus we obtain a, maybe partially def\/ined, map, again denoted by $s:\mathcal{M}(\theta_0,\theta_\infty) \rightarrow
\mathcal{M}(\theta_0',\theta_\infty ')$.
It can be seen from the def\/initions of the elements of ${\rm Aut}({\bf S})$ that the B\"acklund transforma\-tion~$s$ is
a~birational algebraic map.
Using the expression for $s(\tilde{t})$ one obtains a~birational morphism, again denoted by
$s:\tilde{\mathcal{M}}(\theta_0,\theta_\infty) \rightarrow \tilde{\mathcal{M}}(\theta_0',\theta_\infty ')$ which
respects the foliations (i.e., the isomo\-nodromic families).
In particular, $s$ maps `generic' solutions for the parame\-ters~$\theta_0$,~$\theta_\infty$ and the variable~$\tilde{t}$
to `generic' solutions for the parame\-ters~$\theta_0'$,~$\theta_\infty '$ and the variable~$s(\tilde{t})$.
The term `generic' means here the solutions corresponding to irreducible modules and the ones for $(\epsilon_1,\epsilon
_2)$-reducible modules (i.e., Riccati solutions) which are present in both~$\mathcal{M}(\theta_0,\theta_\infty)$ and~$\mathcal{M}(\theta_0',\theta_\infty ')$.

\subsection{Formulas for the B\"acklund transformations}

\subsubsection[$s_1:\tilde{\mathcal{M}}(\theta_0,\theta_\infty)\rightarrow \tilde{\mathcal{M}}(2-\theta_0,-\theta_\infty)$]
{$\boldsymbol{s_1:\tilde{\mathcal{M}}(\theta_0,\theta_\infty)\rightarrow \tilde{\mathcal{M}}(2-\theta_0,-\theta_\infty)}$}

Using the f\/irst charts of the
two spaces and their variables, $s_1$ has the form $(a_{-1},b_{-2},\dots,b_1,c_0,\tilde{t})$ $\mapsto
(a_{-1},b_{-2},\dots,b_1,c_0,\tilde{t}+i\pi)$.
The formula for $s_1$ on the second charts is similar.
The induced map for the ${\rm PIII(D_6)}$ equations is $t\mapsto -t$, $\frac{d}{dt}\mapsto
-\frac{d}{dt}$, $q(\tilde{t})\mapsto q(\tilde{t}+i\pi)$.

\subsubsection[$s_2: \mathcal{M}(\theta_0,\theta_\infty)\rightarrow \mathcal{M}( 1+\theta_0, 1+\theta_\infty)$]
{$\boldsymbol{s_2: \mathcal{M}(\theta_0,\theta_\infty)\rightarrow \mathcal{M}( 1+\theta_0, 1+\theta_\infty)}$}
\label{Section3.4.2}

A point on the f\/irst chart of the
f\/irst space is represented by the dif\/ferential operator $z\frac{d}{dz}+\begin{pmatrix}az^{-1}b\\ z-q-az^{-1}\end{pmatrix}$,
where $a:=a_{-1}$, $q:=-c_0$, $b=b_1z+b_0+b_{-1}z^{-1}+b_{-2}z^{-2}$ and the $b_1,\dots,b_{-2}$
are polynomials in~$t$,~$q$,~$q^{-1}$, using the notation of Section~\ref{Section2.3.1}.
This is transformed by~$s_2$ into the operator $z\frac{d}{dz}+A$ with $A=\begin{pmatrix}az^{-1}+\frac{1}{2}b\\ z-q-az^{-1}+\frac{1}{2}\end{pmatrix}$.
We want to compute an operator $z\frac{d}{dz}+\tilde{A}$
with $\tilde{A}=\begin{pmatrix}\tilde{a}z^{-1}&\tilde{b} \\ z-\tilde{q}&-\tilde{a}z^{-1}\end{pmatrix}$,
$\tilde{b}=\tilde{b}_1z+\tilde{b}_0+ \tilde{b}_{-1}z^{-1}+\tilde{b}_{-2}z^{-2}$ and
$\tilde{b}_1,\dots,\tilde{b}_{-2}$ polynomials in~$t$,~$\tilde{q}$,~$\tilde{q}^{-1}$, representing a~point on the f\/irst chart
of $\mathcal{M}( \frac{1}{2}+\frac{\theta_0}{2}, \frac{1}{2}+\frac{\theta_\infty}{2})$, which is equivalent to
$z\frac{d}{dz}+A$.
Thus we have to solve an equation of the type $\{z\frac{d}{dz}+A\}T=T\{z\frac{d}{dz}+\tilde{A}\}$ with $T\in {\rm GL}(2,\mathbb{C}(z))$.
A local computation shows that~$T$ has the form $T_0+T_{-1}z^{-1}+T_{-2}z^{-2}\neq 0$ with `constant' matrices
$T_0$, $T_{-1}$, $T_{-2}$.
A~Maple computation yields the solution
\begin{gather*}
\tilde{q}=-\frac{tq^2-q\theta_0-t+2a}{q( tq^2+q\theta_\infty -t+2a)},
\qquad
\tilde{a}=\frac{\longg}{ 2q^2( tq^2+q\theta_\infty -t+2a)^2},
\\
\longg=8a^3-4aq^2t^2+8a^2q^2t-qt^2+2aq^4t^2-8a^2t+2at^2-4a^2q+4aqt-q^5t+qt^2\theta_0
\\
\phantom{\longg=}
{}-q^5t^2\theta_0 +q^2t\theta_0^2-4a^2q\theta_0+2aq^2\theta_0+q^4t\theta_0 -q^2t\theta_0-q^5t^2\theta_\infty
-4q^4t\theta_\infty^2+4a^2q\theta_\infty
\\
\phantom{\longg=}
{}+qt^2\theta_\infty-2aq^2\theta_\infty-q^4t\theta_\infty+q^2t^2\theta_\infty+q^3\theta
_0\theta_\infty-4aq^3t\theta_0-4aqt\theta_\infty+q^2t\theta_0\theta_\infty
\\
\phantom{\longg=}
{}-q^4t\theta_0\theta_\infty-2aq^2\theta_0\theta_\infty-4aq^3t+2q^3t.
\end{gather*}

The induced map for solutions of ${\rm PIII(D_6)}$ is obtained from the formula for $\tilde{q}$ and the equality
$q'=\frac{4a-q}{t}$.

{\it Comments on the formulas.}
The term $q$ in the denominator of the formulas is due to our choice of working on the f\/irst charts of the spaces
$\mathcal{M}(\theta_0,\theta_\infty)$ and $\mathcal{M}(\theta_0+1,\theta_\infty +1)$.
This term does not produce singularities for $s_2$.

The denominator $( tq^2+q\theta_\infty -t+2a)$ is due to a~reducible locus of type $(\epsilon_1,\epsilon_2)=(-1,1)$.
More precisely, $tq^2+q\theta_\infty -t+2a=0$ describes the reducible locus of $\mathcal{M}(\theta_0,\theta_\infty)$ if
and only if $\theta_0+\theta_\infty \in 2\mathbb{Z}$ and $\theta_0+\theta_\infty\leq 0$.

{\it A priori}, $s_2$ is not def\/ined on this locus if moreover $\theta_0+\theta_\infty=0$.
One computes that for $\theta_\infty +\theta_0=0$ the formulas reduce to the rational map
$\tilde{q}=-q^{-1}$, $\tilde{a}=\frac{-q+2a}{2q^2}$ and thus $s_2$ is well def\/ined on this locus of type
$(\epsilon_1,\epsilon_2)=(-1,1)$.
In fact, $s_2$ maps this locus to the reducible locus of type $(\epsilon_1,\epsilon_2)=(1,-1)$, which is present in
$\mathcal{M}(\theta_0+1,\theta_\infty +1)$.

\subsubsection[$s_3: \tilde{\mathcal{M}}(\theta_0,\theta_\infty)\rightarrow \tilde{\mathcal{M}}(\theta_0,-\theta_\infty)$]
{$\boldsymbol{s_3: \tilde{\mathcal{M}}(\theta_0,\theta_\infty)\rightarrow \tilde{\mathcal{M}}(\theta_0,-\theta_\infty)}$}

On the f\/irst charts the map is given by:
$(a_{-1},b_{-2},\dots,b_1,c_0,\tilde{t})\mapsto (-ia_{-1},-ib_{-2},b_{-1},ib_0,-b_1$, $-ic_0,\tilde{t}+i\frac{\pi}{2})$ and
similarly on the second charts.
For the ${\rm PIII(D_6)}$ equations, the map is $t\mapsto it$, $\frac{d}{dt}\mapsto -i\frac{d}{dt}$, $q(\tilde{t})\mapsto
-iq(\tilde{t}+i\frac{\pi}{2})$.

\subsubsection[$s_4:\mathcal{M}(\theta_0,\theta_\infty)\rightarrow \mathcal{M}(\theta_\infty,\theta_0)$]
{$\boldsymbol{s_4:\mathcal{M}(\theta_0,\theta_\infty)\rightarrow \mathcal{M}(\theta_\infty,\theta_0)}$}
\label{Section3.4.4}

Consider a~point on the chart $c_1\neq 0$
of $\mathcal{M}(\theta_0,\theta_\infty)$ (see Section~\ref{Section2.3.1}) represented by $z\frac{d}{dz}+A(z)$, where the entries of
$A(z)$ are polynomial in $z$, $z^{-1}$, $t$, $a:=a_{-1}$, $q:=-c_0$, $q^{-1}$.
Now $s_4$ transforms this operator into $z\frac{d}{dz}-A(z^{-1})$.
We suppose that a~transformation $T:=T_{-1}z^{-1}+T_0+T_1z\neq 0$ brings this operator into a~point of the chart
$c_1\neq 0$ of $\mathcal{M}(\theta_\infty,\theta_0)$ represented by an operator $z\frac{d}{dz}+\tilde{A}(z)$, where
the entries of $\tilde{A}(z)$ are polynomials in $z$, $z^{-1}$, $t$, $\tilde{a}=\tilde{a}_{-1}$, $\tilde{q}=-\tilde{c}_0$,
$\tilde{q}^{-1}$.
A~Maple computation shows that there is a~unique solution in terms of $\tilde{q}$ and $\tilde{a}$ of the equation
$\{z\frac{d}{dz}-A(z^{-1})\}T=T\{z\frac{d}{dz}+\tilde{A}(z)\}$, namely
\begin{gather*}
\tilde{q}= \frac{q(-q^2t-\theta_\infty q-t+2a)}{(-q^2t-\theta_0q-t+2a)},
\qquad
\tilde{a} =\frac{\longg}{2(-q^2t-\theta_0q-t+2a)^2},
\\
\longg=4q^2at^2+qt^2\theta_0+q^5t^2\theta_\infty+4q^3at\theta_0-q^5t^2\theta_0-q^4t\theta_\infty\theta_0
-q^2t\theta_0\theta_\infty+4qt\theta_\infty a
\\
\phantom{\longg=}
{}+q^2t\theta_0^2+q^4t\theta_\infty^2-t^2q\theta_\infty
-4q\theta_0a^2-4a^2\theta_\infty q+2a\theta_0q^2\theta_\infty -8a^2q^2t+2aq^4t^2
\\
\phantom{\longg=}
{}+2at^2-8a^2t+8a^3.
\end{gather*}
The induced map for solutions of ${\rm PIII(D_6)}$ is obtained from the formula for $\tilde{q}$ and the equality
$q'=\frac{4a-q}{t}$.

{\it Comments on the formulas.}
The denominator in these formulas is due to a~possible reducible locus of type $(\epsilon_1,\epsilon_2)=(-1,-1)$.
This locus is present in $\mathcal{M}(\theta_0,\theta_\infty)$ if and only if $\theta_0-\theta_\infty \in 2\mathbb{Z}$
and $\theta_0- \theta_\infty \leq 0$.
This locus is not present in $\mathcal{M}(\theta_\infty,\theta_0)$ if moreover $\theta_0-\theta_\infty<0$.
However, in the critical case $\theta_0=\theta_\infty$, $s_4$ turns out to be the identity.

If $\theta_0 -\theta_\infty \in 2\mathbb{Z}$ and $\theta_0-\theta_\infty\geq 2$, then the reducible modules of type
$(1,1)$ are present in $\mathcal{M}(\theta_0,\theta_\infty)$ and not in $\mathcal{M}(\theta_\infty,\theta_0)$.
The corresponding term $2a+tq^2+(2d-1)q+t$ in the denominator of the map $s_4$ does not occur, because $s_4$ maps this
reducible locus to the reducible locus of type $(-1,-1)$ of $\mathcal{M}(\theta_\infty, \theta_0)$.

{\it Using the formulas for $s_1,\dots,s_4$ one can deduce formulas for $B_1$ and $B_2$.
We will however derive these by the direct method used for $s_2$ and $s_4$}.

\subsubsection[The transformation $B_1: \theta_0\mapsto 2+\theta_0$, $\theta_\infty\mapsto\theta_\infty$, $\tilde{t}\mapsto \tilde{t}$]
{The transformation $\boldsymbol{B_1: \theta_0\mapsto 2+\theta_0}$, $\boldsymbol{\theta_\infty\mapsto\theta_\infty}$,
$\boldsymbol{\tilde{t}\mapsto \tilde{t}}$}
\label{Section3.4.5}

We compute the birational map $B_1$ on the open part of the chart $c_1\neq 0$ of
$\mathcal{M}(\theta_0,\theta_\infty)$ where $q=-c_0\neq 0$ (see Section~\ref{Section2.3.1}).
A point is represented by an operator $z\frac{d}{dz}+A$.
The entries of $A$ are polynomials in $q=-c_0$, $q^{-1}$ and $a=a_{-1}$.
Further we suppose that the image under $B_1$ lies in the open part of $\mathcal{M}(2+\theta_0,\theta_\infty)$ given by
$c_1\neq 0$ and has the form $z\frac{d}{dz}+\tilde{A}$, where the entries of $\tilde{A}$ are polynomials in
$\tilde{q}$, $\tilde{q}^{-1}$ and $\tilde{a}=\tilde{a}_{-1}$.
Since the two dif\/ferential operators represent the same dif\/ferential module, there exists a~$T\in {\rm
GL}(2,\mathbb{C}(z))$ such that $z\frac{d}{dz}+\tilde{A}=T^{-1}(z\frac{d}{dz}+A)T$.

Local calculations at $z=0$ and $z=\infty$ predict that $T$ has the form $T=T_{-2}z^{-2}+T_{-1}z^{-1}+T_0$,
$T_{-2}=\begin{pmatrix}0&*\\ 0&0\end{pmatrix}$ and $\det T\in \mathbb{C}^*$.

Further it is assumed that $z\frac{d}{dz}+A$ is not $(\epsilon_1,\epsilon_2)$-reducible for the critical cases
$(\epsilon_1,\epsilon_2)=(-1,-1)$, $\theta_\infty =\theta_0$ and $(\epsilon_1,\epsilon_2)=(-1,1)$, $\theta_\infty
=-\theta_0$, where this reducible locus is present in $\mathcal{M}(\theta_0,\theta_\infty)$ and not in
$\mathcal{M}(2+\theta_0,\theta_\infty)$.
Maple produces the formulas
\begin{gather*}
\tilde{q}=\frac{q(-4a^2+4at-t^2+\theta_0^2q^2+2t\theta_\infty q^3+t^2q^4 )}{(2a-t-\theta_0q+tq^2)(2a-t-\theta_0q-tq^2)},
\\
\tilde{a}=\frac{\longg}{(2a-t-\theta_0q+tq^2)^2(2a-t-\theta_0q-tq^2)^2},
\end{gather*}
where ``long'' means too long for copying.
Substitution of $a=\frac{tq'+q}{4}$ in the formula for $\tilde{q}$ yields the formula for $B_1$ with respect to
solutions.

In the f\/irst critical case $(\epsilon_1,\epsilon_2)=(-1,-1)$, $\theta_\infty =\theta_0$ one obtains
\begin{gather*}
\tilde{q}=-\frac{q(2a-t+\theta_0q+tq^2)}{2a-t-\theta_0q+tq^2},
\qquad
\tilde{a}=\frac{\longg}{(2a-t-\theta_0q+tq^2)^2}.
\end{gather*}
The reducible locus is given by $a=\frac{t}{2}+\frac{\theta_0}{2}q+\frac{t}{2}q^2$.
Thus the above map extends to the reducible locus and produces there the formulas
\begin{gather*}
\tilde{q}=-\frac{q(2tq+2\theta_0+1)}{2tq+1},
\qquad
\tilde{a}=\frac{\longg}{(2tq+1)^2}.
\end{gather*}
For the second critical case $(\epsilon_1,\epsilon_2)=(-1,1)$, $\theta_\infty =-\theta_0$ one f\/inds
\begin{gather*}
\tilde{q}=\frac{q(-2a+t-\theta_0q+tq^2)}{2a-t-\theta_0q-tq^2},
\qquad
\tilde{a}=\frac{\longg}{(2a-t-\theta_0q+tq^2)^2}.
\end{gather*}
The reducible locus is given by $a=\frac{t}{2}+\frac{\theta_0}{2}q-\frac{t}{2}q^2$.
On this locus one has $\tilde{q}=-q+\frac{\theta_0}{t}$.

{\it Comment.}
As in Sections~\ref{Section3.4.2} and~\ref{Section3.4.4},
the map $B_1$ is well def\/ined on the reducible loci because $B_1$ changes the type
$(\epsilon_1,\epsilon_2)$ of the reducible loci.

\subsubsection[The transformation $B_2:\theta_0\mapsto \theta_0$, $\theta_\infty \mapsto 2+\theta_\infty$, $\tilde{t}\mapsto \tilde{t}$]
{The transformation $\boldsymbol{B_2:\theta_0\mapsto \theta_0}$, $\boldsymbol{\theta_\infty \mapsto 2+\theta_\infty}$,
$\boldsymbol{\tilde{t}\mapsto \tilde{t}}$}

Let an object of $\mathcal{M}(\theta_0,\theta_\infty)$ be represented by a~standard operator $z\frac{d}{dz}+A$.
We expect that the transformation $B_2$ yields an object of $\mathcal{M}(\theta_0,2+\theta_\infty)$,
represented by a~standard operator $z\frac{d}{dz}+\tilde{A}$.
Then $z\frac{d}{dz}+\tilde{A}=T^{-1}(z\frac{d}{dz}+A)T$ for a~certain $T\in {\rm GL}(2,\mathbb{C}(z))$.

Local calculations at $z=0$ and $z=\infty$ show that $T$ has the form $T_{-1}z^{-1}+T_0+T_1z$
with $T_{-1}=\begin{pmatrix}0&*\\0&0\end{pmatrix}$, $\det T_1=0$, $\det T\in \mathbb{C}^*$.

The matrix $A$ depends on $q$ and $a:=a_{-1}$ and the matrix $\tilde{A}$ depends on $\tilde{q}$ and
$\tilde{a}:=\tilde{a}_{-1}$.
We have to solve the equation $T(z\frac{d}{dz}+\tilde{A})=(z\frac{d}{dz}+A)T$.

Maple produced the formula
\begin{gather*}
\tilde{q}=-\frac{(2a+t+\theta_\infty q+tq^2)(2a-t+\theta_\infty
q+tq^2)q}{4aq^2t+2qt-t^2-2q^3t-q^2\theta_\infty^2+4a^2-4aq-2qt\theta_0-2q^2\theta_\infty +t^2q^4}.
\end{gather*}
The denominator of $\tilde{a}$ is the square of the denominator of $\tilde{q}$ and the numerator of $\tilde{a}$ is too
large to copy here.
The substitution of $a=\frac{tq'+q}{4}$ in the formula for $\tilde{q}$ yields the $B_2$ map for the solutions of ${\rm PIII(D_6)}$.

The two cases where $B_2$ is, {\it a~priori}, not def\/ined on the reducible locus are:

1.~$(\epsilon_1,\epsilon_2)=(1,1)$ and $\frac{\theta_0}{2}=\frac{\theta_\infty}{2}+1$.
The reducible locus is given by $2a+t+\theta_\infty q+tq^2=0$.
After substitution of $\frac{\theta_0}{2}=\frac{\theta_\infty}{2}+1$, the denominator of $\tilde{q}$ factors as
$(2a+t+\theta_\infty q +tq^2 )(2a-t-\theta_\infty q -2q +tq^2)$. Further $\tilde{q}=-\frac{(2a-t+\theta_\infty q
+tq^2)q}{2a-t-\theta_\infty q -2q+tq^2}$ On the reducible locus one has $\tilde{q}=\frac{tq}{t+\theta_\infty q+q}$ and
$B_2$ is well def\/ined.

2.~$(\epsilon_1,\epsilon_2)=(-1,1)$ and $\frac{\theta_0}{2}=-\frac{\theta_\infty}{2}$.
The reducible locus is given by $2a-t+\theta_\infty q+tq^2=0$.
After substitution of $\frac{\theta_0}{2}=-\frac{\theta_\infty}{2}$, the denominator of $\tilde{q}$ factors as
$(2a-t+\theta_\infty q +tq^2)(2a+t-\theta_\infty q -2q +tq^2)$. Further $\tilde{q}=-\frac{(2a+t+\theta_\infty
q+tq^2)q}{2a+t-\theta_\infty q-2q+tq^2}$.
On the reducible locus one has $\tilde{q}=\frac{tq}{t-\theta_\infty q-q}$ and $B_2$ us well def\/ined.

{\it Comment.}
As in Sections~\ref{Section3.4.2}, \ref{Section3.4.4} and~\ref{Section3.4.5}, the map $B_2$ is well def\/ined because $B_2$ changes the types
$(\epsilon_1,\epsilon_2)$ of the reducible loci.

\subsection*{Acknowledgments}

The authors thank Yousuke Ohyama for his helpful answers to our questions and his remarks concerning the tau-divisor
(see Section~\ref{Section2.3.4}).

%\cite{Gr,Ina06, IIS1, IIS2, IISA,JMU,JM,No,No-Y,OO,O1,O2,O3,O5,STT02,S-Ta,SU01,STe02,Sakai,T}

\pdfbookmark[1]{References}{ref}
\LastPageEnding

\end{document}